\documentclass[12pt,dvips]{article}
\usepackage[T2A]{fontenc}
\usepackage{amssymb,amsmath,amsfonts,amsthm,amscd,latexsym,verbatim,graphics,epsfig,indentfirst,color,stmaryrd}
\usepackage[bookmarks=true,pdfborder=0 0 0,colorlinks=true,linkcolor=blue,dvips=true,final=true]{hyperref}
\usepackage{geometry}
\def\id{\mathrm{id}}

\def\rad{\mathrm{rad}}

\def\rank{\mathrm{rank}}
\def\Span{\mathrm{Span}}

\def\Aut{\mathrm{Aut}}
\def\Ann{\mathrm{Ann}}
\def\Imm{\mathrm{Im}}

\def\tr{\textrm{tr}}

\begin{document}

\begin{center}
{\Large
Rota---Baxter operators of nonzero weight on the matrix algebra of order three}

M. Goncharov, V. Gubarev
\end{center}

\begin{abstract}
We classify all Rota---Baxter operators of nonzero weight
on the matrix algebra of order three over an algebraically closed field of characteristic zero which are not arisen from the decompositions of the entire algebra into a direct vector space sum of two subalgebras.

\end{abstract}

\section{Introduction}

Given an algebra $A$ and a scalar $\lambda\in F$, where $F$ is a~ground field,
a~linear operator $R\colon A\rightarrow A$ is called a Rota---Baxter operator
(RB-operator, for short) on $A$ of weight~$\lambda$ if the following identity
\begin{equation}\label{RB}
R(x)R(y) = R( R(x)y + xR(y) + \lambda xy )
\end{equation}
holds for all $x,y\in A$. Then the algebra $A$ is called a~Rota---Baxter algebra (RB-algebra).

Glen Baxter in 1960~\cite{Baxter} introduced the notion of a Rota---Baxter operator
as formal generalization of integration by parts formula.
Further, F.~Atkinson~\cite{Atkinson}, G.-C.~Rota~\cite{Rota}, P.~Cartier~\cite{Cartier} and others studied such operators on commutative algebras.

At the beginning of the 1980s, the deep connection between solutions of the classical Yang---Baxter equation (named after Rodney Baxter) from mathematical physics and RB-operators on a semisimple finite-dimensional Lie algebra was found by A.A. Belavin and V.G. Drinfel'd~\cite{BelaDrin82} and M.A. Semenov-Tyan-Shanskii~\cite{Semenov83}.
Actually, M.A. Semenov-Tyan-Shanskii had rediscovered the notion of an RB-operator of nonzero weight on a Lie algebra. He called such operators as solutions of modified Yang---Baxter equation.

From the 2000s, the active study of Rota---Baxter operators on associative algebras has  begun, see the monograph~\cite{GuoMonograph}.
To the moment, many applications and connections of Rota---Baxter operators
with symmetric polynomials, shuffle algebra, Loday algebras, etc. were found \cite{Atkinson,Shuffle,Aguiar00,BurdeGubarev,GuoMonograph}.

One of the interesting direction in the study of Rota---Baxter operators is a problem of classification of RB-operators on a given algebra. RB-operators were classified on $\mathrm{sl}_2(\mathbb{C})$
\cite{Kolesnikov,KonovDissert,sl2,sl2-0}, on $M_2(\mathbb{C})$ \cite{BGP,Mat2},
on $\mathrm{sl}_3(\mathbb{C})$ \cite{KonovDissert} and on other algebras~\cite{AnBai,BGP,3-Lie}.
In 2013, V.V.~Sokolov described all skew-symmetric RB-operators of nonzero weight on $M_3(\mathbb{C})$~\cite{Sokolov}; up to conjugation with automorphisms and transpose he got 8 series.

In~\cite{Unital}, some general properties of RB-operators were stated. In particular, every RB-operator $R$ of weight zero on $M_n(F)$ over a field $F$ of characteritic zero is nilpotent and $R^{2n-1} = 0$.
Given an RB-operator~$R$ of nonzero weight on $M_n(\mathbb{C})$, we may assume that $R(1)$ is diagonal (up to conjugation with automorphisms of $M_n(\mathbb{C})$).
The last result gives a powerful tool for the study of RB-operators of nonzero weight on the matrix algebra. As a corollary, it was proved in~\cite{Unital} that every RB-operator $R$ of nonzero weight on $M_3(\mathbb{C})$ is diagonal, it means that the subalgebra of diagonal matrices in $M_3(\mathbb{C})$ is preserved by $R$ (up to conjugation with automorphisms).
In the current work we classify all RB-operators on $M_3(\mathbb{C})$ which are not projections of $M_3(\mathbb{C})$ onto a~subalgebra parallel to another one. The study of such projections is a part of the research area of the decompositions of an algebra into a sum of two subalgebras, see, e.g.~\cite{BurdeGubarev}.

Let us remark that for the algebra $M_n(\mathbb{C})$, solutions of the associative Yang---Baxter equation~\cite{Zhelyabin,Aguiar} are in one-to-one correspondence with RB-operators of weight zero~\cite{Unital}. The same bijection holds~\cite{AYBE-ext} for solutions of the weighted associative Yang---Baxter equation~\cite[p.~113]{FardThesis} and RB-operators of nonzero weight on $M_n(\mathbb{C})$.

Let us give a brief outline of the work.
In~\S2, we give some required preliminaries in Rota---Baxter operators.
In~\S3, we consider some results about RB-operators of nonzero weight on the matrix algebra.
In Corollary~1, we clarify the classification of RB-operators on $M_2(\mathbb{C})$ from~\cite{BGP} .

As we said above, an RB-operator of nonzero weight on $M_3(\mathbb{C})$
preserves the subalgebra of diagonal matrices. In \S4, we apply the descriptions of RB-operators of nonzero weight on $F\oplus F\oplus F$ from~\cite{SumOfFields}.

In~\S5, we prove the main result. In Theorem~3, we classify all RB-operators of nonzero weight on $M_3(\mathbb{C})$, we get 36 cases up to all natural actions.
Given an RB-operator $R$ of weight~1 on the subalgebra $D_3$ of diagonal matrices,
we may extend it on the entire algebra $M_3(\mathbb{C})$ as follows:
we put $R(U_3) = 0$ and $(R+\id)(L_3) = 0$, where $U_3$ and $L_3$ denote the subalgebra of upper and lower triangular matrices respectively. Thus, we get 20 cases~in~Theorem~3.
In the same way we may extend a given RB-operator $R$ of $F\oplus M_2(\mathbb{C})$
onto the entire $M_3(\mathbb{C})$; such construction gives another 12 cases in Theorem~3.
Finally, we also have four ``exceptional'' cases in Theorem~3.
In Corollary~2, we check that all obtained RB-operators lie in pairwise distinct orbits.

\section{Preliminaries}

{\bf Statement 1} \cite{GuoMonograph}.
Given an RB-operator $P$ of weight $\lambda$,

a) the operator $-P-\lambda\id$ is again an RB-operator of weight $\lambda$,

b) the operator $\lambda^{-1}P$ is an RB-operator of weight 1, provided $\lambda\neq0$.

Given an algebra $A$, let us define a map $\phi$ on the set of all RB-operators on $A$
as $\phi(P)=P'=-P-\lambda(P)\id$, where $\lambda(P)$ denotes the weight of an RB-operator~$P$. Note that $\phi^2$ coincides with the identity map.

{\bf Statement 2} \cite{BGP}.
Given an algebra $A$, an RB-operator $P$ of weight $\lambda$ on $A$,
and $\psi\in\Aut(A)$, the operator $P^{(\psi)} = \psi^{-1}P\psi$
is an RB-operator of weight $\lambda$ on $A$.

{\bf Statement 3} \cite{GuoMonograph}.
Let an algebra $A$ to split as a vector space
into a direct sum of two subalgebras $A_1$ and $A_2$.
An operator $P$ defined as
\begin{equation}\label{Split}
P(a_1 + a_2) = -\lambda a_2,\quad a_1\in A_1,\ a_2\in A_2,
\end{equation}
is an RB-operator of weight $\lambda$ on $A$.

Let us call an RB-operator from Statement~3 as
{\it splitting} RB-operator with subalgebras $A_1,A_2$.
There is a~bijection between the set of all splitting RB-operators on
an algebra~$A$ and all decompositions of~$A$ into a direct sum of two subalgebras $A_1,A_2$.

Note that if $P$ is a splitting RB-operator on $A$ of weight
$\lambda$ with subalgebras $A_1,A_2$, then
$\phi(P)$ is the splitting RB-operator of weight $\lambda$
with the same subalgebras $A_1, A_2$ (just another projection).

{\bf Statement 4} \cite{BGP}.
Let $A$ be a unital algebra, and let $P$ be an RB-operator of nonzero weight~$\lambda$ on $A$.

a) If $P(1)\in F$, then $P$ is splitting;

b) If $P(P(x) + \lambda x) = 0$ for all $x\in A$, then $P$ is splitting.

Given a unital algebra~$A$, we call an RB-operator $R$ on~$A$ of nonzero weight
as {\it inner-splitting} if $R(1)\in F$.

{\bf Statement 5} \cite{Guil}.
Let an algebra $A$ be a direct sum of subalgebras
$A_-,A_0,A_+$, and $A_\pm$ are $A_0$-modules.
If $R_0$ is an RB-operator of weight~$\lambda$ on $A_0$,
then an operator $P$ defined as follows
\begin{equation}\label{RB:SubAlg2}
P(a_-+a_0+a_+) = R_0(a_0) - \lambda a_+,\quad
a_{\pm}\in A_{\pm},\ a_0\in A_0,
\end{equation}
is an RB-operator of weight $\lambda$ on $A$.

Let us call an RB-operator of nonzero weight
defined by \eqref{RB:SubAlg2} as {\it triangular-splitting}
provided that at least one of $A_-,A_+$ is nonzero.
If $A_0 = (0)$, then $P$ is splitting RB-operator on $A$.
If $A_0$ has trivial (zero) product, then any linear map on $A_0$ is suitable as $R_0$.
Note that if $P$ is a triangular-splitting RB-operator on an algebra $A$
with subalgebras $A_\pm,A_0$, then the operator $\phi(P)$ is the
triangular-splitting RB-operator with the same subalgebras.

{\bf Statement 6}~\cite{Unital}.
Let an algebra $A$ be equal a direct sum of two ideals $A_1$ and $A_2$
and $R$ be an RB-operator of weight~$\lambda$ on $A$.
Then $\mathrm{Pr}_i R$ is the RB-operator of weight~$\lambda$ on $A_i$,
$i=1,2$. Here $\mathrm{Pr}_i$ denotes the projection from $A$ onto $A_i$.

\section{RB-operators on $M_n(F)$}

Over a field~$F$, denote the algebra of all upper and lower triangular matrices of order~$n$ by $U_n$ and $L_n$ respectively.

{\bf Example 1}~\cite{BGP}.
Decomposing $M_n(F) = L_n\oplus D_n\oplus U_n$
(as vector spaces) and given an RB-operator defined on $D_n$, we have a~triangular-splitting RB-operator defined with $A_- = U_n$, $A_+ = L_n$, $A_0 = D_n$.

{\bf Example 2}.
Let $1\leq k<n$, define
\begin{gather*}
M_1 = \Span\{e_{ij}\mid 1\leq i,j\leq k\}\oplus \Span\{e_{ij}\mid k+1\leq i,j\leq n\},\\
M_2 = \Span\{e_{ij}\mid 1\leq i\leq k;\,k+1\leq j\leq n\}, \\
M_3 = \Span\{e_{ij}\mid k+1\leq i\leq n;\,1\leq j\leq k\}.
\end{gather*}
Decomposing $M_n(F) = M_1\oplus M_2\oplus M_3$ (as vector spaces) and given an RB-operator defined on $M_1$,
we have a~triangular-splitting RB-operator on $M_n(F)$ defined with
$A_+ = M_2$, $A_- = M_3$, $A_0 = M_1$.

{\bf Example 3}.
Let an algebra~$A$ be equal a direct vector space sum $A_- \oplus A_0 \oplus A_+$
of its three subalgebras satisfying the conditions $A_0 A_-,A_-A_0\subset A_+$.
Suppose that~$P$ is an RB-operator of weight~1 on $A_0$.
Then a~linear operator $R$ defined on~$A$ as
$R(a_- + a_0 + a_+) = P(a_0) - a_+$
is an RB-operator of weight~1 if $R(A_0)A_-,A_-R(A_0)\subset A_-$.

{\bf Lemma 1}.
Let $R$ be an RB-operator of weight~1 on $M_n(F)$ such that $R(1)$ is a~diagonal matrix,
\begin{equation}\label{R(1)diagonal}
\begin{gathered}
R(1) = \sum\limits_{i=1}^s \lambda_i
(e_{k_1+\ldots+k_{i-1}+1,k_1+\ldots+k_{i-1}+1}+\ldots+e_{k_1+\ldots+k_i,k_1+\ldots+k_i}),\\
k_0 = 0,\quad k_1,\ldots,k_s\geq 1,\quad \sum\limits_{i=1}^s k_i = n.
\end{gathered}
\end{equation}
Suppose that either $s$-tuple
$(\lambda_1,\ldots,\lambda_s)$ or $(\lambda_s,\ldots,\lambda_1)$
equals $(-f,-f+1,\ldots,0,\ldots,g)$ for $f,g\geq 0$.
Then the space
$$
V_t = \Span\{e_{ij}\mid k_1+\ldots+k_{p-1}<i\leq k_1+\ldots+k_p,
                      \,k_1+\ldots+k_{q-1}<j\leq k_1+\ldots+k_q,\,p-q=t\}
$$
is $R$-invariant for all $t$.
Moreover, $R$ is splitting on both $V_{f+g}$ and $V_{-f-g}$.

{\sc Proof}.
Write down the following equalities,
\begin{gather}
R(1)R(x) = R^2(x) + R(x) + R(R(1)x), \label{R(1)R(x)} \\
R(x)R(1) = R^2(x) + R(x) + R(xR(1)). \label{R(x)R(1)}
\end{gather}
From~\eqref{R(1)R(x)} and~\eqref{R(x)R(1)}, we deduce that
\begin{equation}\label{[R(1),R(x)]}
[R(1),R(x)] = R([R(1),x]).
\end{equation}

Let $x$ be a matrix unity from $V_t$, it means $x = e_{yz}$
for some $k_1+\ldots+k_{p-1}<y\leq k_1+\ldots+k_p$ and
$k_1+\ldots+k_{q-1}<z\leq k_1+\ldots+k_q$ with $p - q = t$.
For $i,j$ such that
$k_1+\ldots+k_{u-1}<i\leq k_1+\ldots+k_u$ and
$k_1+\ldots+k_{v-1}<j\leq k_1+\ldots+k_v$ and
for $A = (a_{cd})_{c,d=1}^n = R(x)$,
we get the equalities
$a_{ij}(\lambda_u-\lambda_v) = a_{ij}(\lambda_p-\lambda_q)$ by~\eqref{[R(1),R(x)]}.
When $u - v \neq t$, we get $a_{ij} = 0$.

Consider $V_{f+g}$, it is exactly the block
$\Span\{e_{ij}\mid 1\leq i\leq k_1,\,k_1+\ldots+k_{s-1}<j\leq k_1+\ldots+k_s\}$.
Putting a~matrix unity $e_{yz}$ from $V_{f+g}$ into~\eqref{R(1)R(x)},
we get $R^2(e_{yz}) + R(e_{yz}) = 0$. By Statement~4b), $R$~is splitting on $V_{f+g}$.
Analogously, $R$~is splitting on $V_{-f-g}$.
\hfill $\square$

{\bf Remark 1}.
Actually, in~\cite{Unital} it was proved that $V_0$ is $R$-invariant.
We have extended this statement in Lemma~1 for all~$t$.

{\bf Theorem 1}~\cite{BGP,Unital}.
Let $F$ be an algebraically closed field of characteristic zero.
Every nontrivial RB-operator of weight~1 on $M_2(F)$ up to conjugation with an automorphism of $M_2(F)$ or transpose, up to~$\phi$ equals one of the following cases:

(a) $R\begin{pmatrix}
x_{11} & x_{12} \\
x_{21} & x_{22} \\
\end{pmatrix}
 = \begin{pmatrix}
0 & - x_{12} \\
0 & x_{11}
\end{pmatrix}$,

(b) $R\begin{pmatrix}
x_{11} & x_{12} \\
x_{21} & x_{22} \\
\end{pmatrix}
 = \begin{pmatrix}
-x_{11} & - x_{12} \\
0 & 0
\end{pmatrix}$,

(c) $R\begin{pmatrix}
x_{11} & x_{12} \\
x_{21} & x_{22} \\
\end{pmatrix}
 = -x_{21}\begin{pmatrix}
\alpha & \alpha\gamma \\
1 & \gamma
\end{pmatrix}$,\quad $\alpha,\gamma \in F$,

(d) $R\begin{pmatrix}
x_{11} & x_{12} \\
x_{21} & x_{22} \\
\end{pmatrix}
 = \begin{pmatrix}
\alpha x_{12} & - x_{12}\\
-x_{21} & (1/\alpha)x_{21}
\end{pmatrix}$,\quad $\alpha\in F$,

(e) $R\begin{pmatrix}
x_{11} & x_{12} \\
x_{21} & x_{22} \\
\end{pmatrix}
 = \begin{pmatrix}
x_{22}-x_{11} + \beta x_{21} & (-\beta^2/4) x_{21} - x_{12}\\
0 & 0
\end{pmatrix}$,\quad $\beta\in F$.

Let us refine the classification of RB-operators on $M_2(F)$
from Theorem~1 in the following way.

{\bf Corollary 1}.
Let $F$ be an algebraically closed field of characteristic zero.
Every nontrivial RB-operator on $M_2(F)$ of weight~1,
up to conjugation with an automorphism of $M_2(F)$ or transpose and up to the action of $\phi$ from Statement 1, equals one of the following cases:

(M1) $R\begin{pmatrix}
x_{11} & x_{12} \\
x_{21} & x_{22} \\
\end{pmatrix}
 = \begin{pmatrix}
0 & - x_{12} \\
0 & x_{11}
\end{pmatrix}$,

(M2) $R\begin{pmatrix}
x_{11} & x_{12} \\
x_{21} & x_{22} \\
\end{pmatrix}
 = \begin{pmatrix}
-x_{11} & - x_{12} \\
0 & 0
\end{pmatrix}$,

(M3) $R\begin{pmatrix}
x_{11} & x_{12} \\
x_{21} & x_{22} \\
\end{pmatrix}
 = \begin{pmatrix}
0 & 0 \\
-x_{21} & 0
\end{pmatrix}$,

(M4) $R\begin{pmatrix}
x_{11} & x_{12} \\
x_{21} & x_{22} \\
\end{pmatrix}
 = \begin{pmatrix}
x_{21} & 0 \\
-x_{21} & 0
\end{pmatrix}$,

(M5) $R\begin{pmatrix}
x_{11} & x_{12} \\
x_{21} & x_{22} \\
\end{pmatrix}
 = \begin{pmatrix}
x_{12} & - x_{12}\\
-x_{21} &  x_{21}
\end{pmatrix}$,

(M6) $R\begin{pmatrix}
x_{11} & x_{12} \\
x_{21} & x_{22} \\
\end{pmatrix}
 = \begin{pmatrix}
x_{22}-x_{11} & - x_{12}\\
0 & 0
\end{pmatrix}$.

\noindent Moreover, these 6 cases lie in different orbits
of the set of RB-operators of weight~1 on $M_2(F)$
under the action of $\phi$ and conjugation with $\Aut(M_2(F))$ or transpose.

{\sc Proof}.
We are applying the classification from Theorem~1.
The RB-operator $R$ from (c) when $\alpha+\gamma=0$
is conjugate to the RB-operator $P$ from (M3) with the help of the automorphism
\begin{equation}\label{AutForM2}
\begin{gathered}
\psi(e_{11}) = e_{11} - \alpha e_{12},\quad
\psi(e_{12}) = e_{12},\\
\psi(e_{22}) = e_{22} + \alpha e_{12},\quad
\psi(e_{21}) = \alpha(e_{11} - e_{22}) - \alpha^2 e_{12} + e_{21},
\end{gathered}
\end{equation}
it means that $\psi^{-1}P\psi = R$.
When $\alpha+\gamma\neq0$, $R$ is conjugate to (M4) under the action of the automorphism
\begin{gather*}
\chi(e_{11}) = e_{11} + \gamma e_{12},\quad
\chi(e_{12}) = -(\alpha+\gamma)e_{12},\quad
\chi(e_{22}) = e_{22} - \gamma e_{12},\\
\chi(e_{21}) = \frac{1}{\alpha+\gamma}( \gamma(e_{11} - e_{22}) + \gamma^2 e_{12} - e_{21} ).
\end{gather*}

Further, the RB-operator $R$ from (d) is conjugate to (M5) with the automorphism
$$
\xi(e_{11}) = e_{11},\quad
\xi(e_{12}) = (1/\alpha) e_{12},\quad
\xi(e_{22}) = e_{22},\quad
\xi(e_{21}) = \alpha e_{21}.
$$

The RB-operator $R$ from (e) is conjugate to (M6)
with the automorphism~$\psi$~\eqref{AutForM2} defined with $\alpha = \beta/2$.

Now, let us clarify that all 6 cases lie in different orbits.
We note that given an RB-operator $R$,
conjugation with an automorphism or transpose does not change
the algebraic properties of $\ker(R)$ and $\ker(R+\id)$.
The case (M1) is unique non-splitting RB-operator from the list.
The case (M2) is unique splitting RB-operator with $R(1)\not \in F$.
The cases (M3) and (M4) satisfy the condition $\dim(\ker(R)) = 3$
in contrast to (M5) and~(M6). We may distinguish the cases (M3) and (M4)
as follows: the algebra $\ker(R+\id)$ is nilpotent in (M3) but not in (M4).
Finally, in the case (M5), in contrast to the case (M6), we have one of the kernels
isomorphic to $F\oplus F$.
\hfill $\square$

{\bf Remark 2}.
It is easy to show that the RB-operator (M4) is conjugate to the RB-operator
(M4${}^\prime$) $R\begin{pmatrix}
x_{11} & x_{12} \\
x_{21} & x_{22} \\
\end{pmatrix}
 = \begin{pmatrix}
0 & 0 \\
-x_{21} & x_{21}
\end{pmatrix}$.

Let us call an RB-operator $R$ on $M_n(F)$ {\it diagonal}, if
$R^{(\psi)}(D_n)\subseteq D_n$ for some $\psi\in\Aut(M_n(F))$.
By Theorem~1, all RB-operators of nonzero weight on $M_2(F)$ are diagonal.
In~\cite{Unital}, the same result was stated for the matrix algebra of order three.

{\bf Theorem~2}~\cite{Unital}.
Let $F$ be an algebraically closed field of characteristic zero.
All RB-operators of nonzero weight on $M_3(F)$ are diagonal.

In advance, we will apply the following automorphism $\Phi_{ab}$ of $M_3(F)$,
where $a,b\in\{1,2,3\}$ and $a\neq b$.
Let $c$ be such that $\{a,b,c\} = \{1,2,3\}$.
The $\Phi_{ab}$ acts on matrix unities as $\Phi_{ab}(e_{ij}) = e_{i'j'}$,
where $a' = b$, $b' = a$ and $c' = c$.
Moreover, we define an automorphism $\Phi_{abc}$ of $M_3(F)$ for $\{a,b,c\} = \{1,2,3\}$
acting on the indices of matrix unities as the cycle $(abc)$.

\section{RB-operators on $F^3$}

In~\cite{SumOfFields}, RB-operators of nonzero weight on a sum
of fields were studied. In particular, it was proved that

{\bf Statement~7}~\cite{SumOfFields}.
Let $R$ be a not inner-splitting RB-operator of weight~1 on the sum of fields
$F^3 = Ff_1\oplus Ff_2\oplus Ff_3$,
then up to permutation of coordinates and action of $\phi$, we have nine cases:
\begin{align*}
1&) & R(f_1) &=  f_2 + f_3, & R(f_2) &= f_3,        & R(f_3) &= 0, \\
2&) & R(f_1) &=  f_2 + f_3, & R(f_2) &= f_3,        & R(f_3) &= -f_3, \\
3&) & R(f_1) &=  f_2 + f_3, & R(f_2) &= -(f_2+f_3), & R(f_3) &= -f_3, \\
4&) & R(f_1) &=  f_2 + f_3, & R(f_2) &= -f_2,       & R(f_3) &= 0, \\
5&) & R(f_1) &=  f_2 + f_3, & R(f_2) &= 0,          & R(f_3) &= 0, \\
6&) & R(f_1) &=  f_2,       & R(f_2) &= 0,          & R(f_3) &= 0, \\
7&) & R(f_1) &=  f_2,       & R(f_2) &= 0,          & R(f_3) &= -f_3, \\
8&) & R(f_1) &= -f_1,       & R(f_2) &= 0,          & R(f_3) &= 0,\\
9&) & R(f_1) &=  f_2,       & R(f_2) &= -f_2,       & R(f_3) &= -f_3.
\end{align*}

Let us also write down the action of the RB-operator $R' = -(R+\id)$ in all nine cases,
\begin{align*}
1&) & R'(f_1) &= -f_1 - f_2 - f_3, & R'(f_2) &= -f_2 - f_3, & R'(f_3) &= -f_3, \\
2&) & R'(f_1) &= -f_1 - f_2 - f_3, & R'(f_2) &= -f_2 - f_3, & R'(f_3) &= 0, \\
3&) & R'(f_1) &= -f_1 - f_2 - f_3, & R'(f_2) &= f_3,        & R'(f_3) &= 0, \\
4&) & R'(f_1) &= -f_1 - f_2 - f_3, & R'(f_2) &= 0,          & R'(f_3) &= -f_3, \\
5&) & R'(f_1) &= -f_1 - f_2 - f_3, & R'(f_2) &= -f_2,       & R'(f_3) &= -f_3, \\
6&) & R'(f_1) &= -f_1 - f_2,       & R'(f_2) &= -f_2,       & R'(f_3) &= -f_3, \\
\end{align*}
\begin{align*}
7&) & R'(f_1) &= -f_1 - f_2,       & R'(f_2) &= -f_2,       & R'(f_3) &= 0, \\
8&) & R'(f_1) &= 0,                & R'(f_2) &= -f_2,       & R'(f_3) &= -f_3,\\
9&) & R'(f_1) &= -f_1 - f_2,       & R'(f_2) &= 0,          & R'(f_3) &= 0.
\end{align*}

Now, we state some required lemmas about RB-operators of weight~1 on $M_3(F)$.
By~$E$ we denote the unit of $M_3(F)$.

{\bf Lemma 2}.
If there exists $i$ such that $R(e_{ii})\in \{0,-E\}$, then
$$
e_{ii}\Imm(R+\id),\Imm(R+\id)e_{ii}\subset \ker(R).
$$

{\sc Proof}.
Suppose that $R(e_{ii}) = 0$ or $R(e_{ii}) = -E$. For all $x$ we have
\begin{equation}\label{ImR'InKerR}
R( e_{ii}(R(x)+x) )
 = R(e_{ii})R(x) - R(R(e_{ii})x)
 = 0.
\end{equation}
We deal analogously with $\Imm(R+\id)e_{ii}$. \hfill $\square$

We call a subspace $V$ of $M_3(F)$ {\it homogeneous}
if for any $a = \sum\limits_{i,j=1}^3 \alpha_{ij}e_{ij}\in V$
we have $\alpha_{ij}e_{ij}\in V$ for all $i,j=1,2,3$.

{\bf Lemma 3}.
If there exist $i\neq j$ such that $R(e_{ii}),R(e_{jj})\in \{0,-E\}$,
then $\ker(R)$ is homogeneous.

{\sc Proof}.
Since $\ker(R)\subset \Imm(R+\id)$, by Lemma~2
we get that
\begin{equation}\label{Homogeneous-Cond}
e_{ii}\ker(R),e_{jj}\ker(R),\ker(R)e_{ii},\ker(R)e_{jj}\subset \ker(R).
\end{equation}
Let $\{i,j,k\} = \{1,2,3\}$, so,
$e_{kk}\ker(R) = (1-e_{ii}-e_{jj})\ker(R)\subset \ker(R)$
and $\ker(R)e_{kk}\subset \ker(R)$.
If $a\in \ker(R)$, then $e_{ss} a e_{tt}\in \ker(R)$ for all $s,t=1,2,3$. \hfill $\square$

{\bf Lemma 4}.
If there exist $i\neq j$ such that $R(e_{ii})\in \{0,-E\}$
and $R(e_{jj}) = e_{ii}$ or $R(e_{jj}) = e_{ii} - 1$,
then $\ker(R)$ is homogeneous.

{\sc Proof}.
Suppose that $R(e_{jj}) = e_{ii}$ or $R(e_{jj}) = e_{ii} - 1$.
Applying~\eqref{ImR'InKerR}, we get
\begin{multline*}
R( e_{jj}(R(x)+x) )
 = R(e_{jj})R(x) - R(R(e_{jj})x)
 = e_{ii}R(x) - R(e_{ii} x) \\
 = e_{ii}R(x) + R(e_{ii}R(x)),
\end{multline*}
so $e_{jj}\ker(R),\ker(R)e_{jj}\subset \ker(R)$.
By Lemma~2, $e_{ii}\ker(R),\ker(R)e_{ii}\subset \ker(R)$.
So, we get~\eqref{Homogeneous-Cond} and thus $\ker R$ is homogeneous.
\hfill $\square$

Given an RB-operator of weight~1 on $M_3(F)$, by Theorem~2
we may assume that $D_3(F)$ is $R$-invariant and the action of~$R$ on $D_3(F)$
is one of Cases~1--9 listed in Statement~7.

Let us apply Lemma~3 (for Cases 2--5, 8) and Lemma~4 (for 1, 6, 7) to get the following data about RB-operators of weight~1 on $M_3(F)$ with prescribed action
on~$D_3(F)$ (see the table on the next page).

{\bf Lemma 5}.
Let $R$ be an RB-operator on $M_3(F)$ of weight~1.
Suppose that $R'(e_{ii}) = -E$ for some $i\in\{1,2,3\}$.
Then $\Imm(R)$ has zero projection on $e_{ii}$.

{\sc Proof}.
By Lemma~2, $e_{ii}\Imm(R)e_{ii}\subset \ker(R')$.
Since $R'(e_{ii})\neq 0$, the space $\Imm(R)$ has zero projection on $e_{ii}$. \hfill $\square$

\begin{center}
\begin{tabular}{c|c|c}
Case & $\ker(R)$-homogeneity & $\ker(R')$-homogeneity \\
\hline
1 & + & + \\
2 &   & + \\
3 &   & + \\
4 &   & + \\
5 & + &  \\
6 & + &  \\
7 & + & + \\
8 & + &  \\
9 &   & + \\
\end{tabular}
\end{center}

\section{Main result}

Call an RB-operator of weight~1 on $M_3(F)$ defined by Example~1 as {\it primitive} one.

Let $R$ be a primitive RB-operator on $M_3(F)$ and its action on $D_3(F)$
corresponds to Case $X$, where $X\in\{1,2,3,4,5,6,7,8,9\}$ (see Statement~7).
We denote it as case~$Xa$). By case $Xb$) ($Xc$)), we denote
a primitive RB-operator which action on $D_3(F)$ is conjugate with the automorphism $\Phi_{12}$ ($\Phi_{23}$) (restricted on~$D_3(F)$) to Case $X$.

{\bf Theorem 3}.
Let $F$ be an algebraically closed field of characteristic zero.
Every non-splitting RB-operator of weight~1 on $M_3(F)$,
up to conjugation with an automorphism of $M_3(F)$ or transpose and up to the action of $\phi$

A) is defined by Example~1 with all cases 1a)--7a), 1b)--7b), 1c)--7c) except 5c);

B) or is one of the following ones

\noindent 1-I) $R(e_{11}) = e_{22}+e_{33}$, $R(e_{22}) = e_{33}$;
$e_{31},e_{32},e_{33}\in \ker(R)$,
$e_{12},e_{13}\in\ker(R')$,
$R(e_{21})=-e_{32}$, $R(e_{23})=e_{12}-e_{23}$;

\noindent 6-I)
$e_{22},e_{23},e_{32},e_{33}\in\ker(R)$,
$R(e_{11}) = e_{22}$,
$R(e_{12})=-e_{12}-e_{32}$,
$R(e_{21})=-e_{21}+e_{23}$,
$R(e_{13})=e_{11}+e_{22}-e_{13}$,
$R(e_{31})=-e_{31}+e_{33}$;

\noindent 6-II)
$e_{12},e_{13},e_{22},e_{23},e_{32},e_{33}\in\ker(R)$,
$R(e_{11}) = e_{22}$, $R(e_{21})=-e_{21}$,
$R(e_{31})=-e_{31}+e_{11}+e_{22}$;

\noindent 6-III)
$e_{12},e_{13},e_{22},e_{23},e_{32},e_{33}\in\ker(R)$,
$R(e_{11}) = e_{22}$,
$R(e_{21})=-e_{21}-e_{23}$,
$R(e_{31})=-e_{31}+e_{11}+e_{22}$;

C) or is defined by Example~2 for $k=1$, i.e.,
$e_{12},e_{13}\in\ker(R')$,
$e_{21},e_{31}\in\ker(R)$,

\noindent2-I) $R(e_{11}) = e_{22}+e_{33}$, $R(e_{22}) = e_{33}$;
$e_{23},e_{33}\in \ker(R')$,
$R(e_{32}) = e_{33}$;

\noindent2-II) $R(e_{11}) = e_{22}+e_{33}$, $R(e_{22}) = e_{33}$;
$e_{32},e_{33}\in \ker(R')$,
$R(e_{23}) = e_{33}$;

\noindent3-I) $e_{11},e_{23}\in \ker(R')$, $-R(e_{22}) = R(e_{33}) = e_{11}+e_{22}$;
$R(e_{32}) = e_{11}+e_{22}$;

\noindent3-II) $e_{11},e_{32}\in \ker(R')$, $-R(e_{22}) = R(e_{33}) = e_{11}+e_{22}$;
$R(e_{23}) = e_{11}+e_{22}$;

\noindent4-I) $R(e_{22}) = e_{11}+e_{33}$;
$e_{33},e_{32}\in \ker(R')$,
$e_{11}\in \ker(R)$,
$R(e_{23}) = e_{33}$;

\noindent4-II) $R(e_{22}) = e_{11}+e_{33}$;
$e_{33},e_{23}\in \ker(R')$,
$e_{11}\in \ker(R)$,
$R(e_{32}) = e_{33}$;

\noindent5-I) $R(e_{11}) = e_{22}+e_{33}$;
$e_{22},e_{23},e_{32},e_{33}\in \ker(R)$;

\noindent5-II)
$R(e_{11}) = e_{22}+e_{33}$;
$e_{22},e_{23},e_{33}\in \ker(R)$,
$R(e_{32}) =  e_{33} - e_{32}$;

\noindent6-IV) $e_{11},e_{33},e_{32}\in\ker(R)$, $R(e_{22}) = e_{11}$,
$R(e_{23})=-e_{23} + e_{11}+e_{22}$;

\noindent6-V) $e_{11},e_{33},e_{23}\in\ker(R)$, $R(e_{22}) = e_{11}$,
$R(e_{32})=-e_{32} + e_{11}+e_{22}$;

\noindent 6-VI)
$e_{22},e_{33},e_{23},e_{32}\in\ker(R)$,
$R(e_{11}) = e_{22}$;

\noindent8-I) $e_{11}\in \ker(R')$,
$e_{22},e_{23},e_{33}\in \ker(R)$,
$R(e_{32}) = e_{11} + e_{22} - e_{32}$.

{\sc Proof}.
Let $R$ be an RB-operator of weight~1 on $M_3(F)$. By Theorem~2, we assume that $R$ is diagonal and we have one of Cases 1)--9) from Statement~7 for the action of~$R$ on $D_3(F)$. To prove Theorem, we consider all of them case-by-case.

Let us show that all primitive non-splitting RB-operators on $M_3(F)$
are described in A).
First, RB-operators defined by Example~1 for Cases 8) and 9) are splitting.
Further, let~$R$ be an RB-operator on $M_3(F)$ defined by Example~1 and its action on $D_3(F)$
is conjugate with the automorphism $\Phi_{13}$ (restricted on~$D_3(F)$) to Case $X$, $X\in\{1,2,3,4,5,6,7\}$.
Up to transpose, $R$ is conjugate with the automorphism $\Phi_{13}$ of~$M_3(F)$
to Case $X$. So, the action of $S_3 = \Aut(D_3(F))$ gives only three different cases $Xa$, $Xb$, and $Xc$ for each~$X$.
Finally, the subcase 5c) coincides with the subcase 5b).

Now we consider one by one possible actions  of $R$ on $D_3(F)$ from Statement 7.

\underline{Case 1)} $R(e_{11}) = e_{22} + e_{33}$, $R(e_{22}) = e_{33}$, $R(e_{33}) = 0$.

In this case $R(1) = e_{22} + 2e_{33}$.
Now we want to use Lemma~1 to specify the information about~$R$. Let us illustrate the statement of Lemma~1 in this case in details. We will need the equation \eqref{[R(1),R(x)]}. Let $a\in M_3(F)$ and let $R(a) = (\alpha_{ij})$. Then
$$
[R(1),R(a)]=\left (
\begin{array}{ccc}
0  & -\alpha_{12} & -2\alpha_{13}    \\
\alpha_{21}  & 0  & -\alpha_{23}  \\
2\alpha_{31}  & \alpha_{32} & 0
\end{array} \right).
$$

Now we can one by one put $e_{ij}$ instead of $a$. First take
$a = e_{12}$. Since
$[R(1),e_{12}] = [e_{22}+2e_{33},e_{12}] = -e_{12}$,
we obtain
$$
\left (
\begin{array}{ccc}
0  & -\alpha_{12} & -2\alpha_{13}    \\
\alpha_{21}  & 0  & -\alpha_{23}  \\
2\alpha_{31}  & \alpha_{32} & 0
\end{array} \right)=-\left (
\begin{array}{ccc}
\alpha_{11}  & \alpha_{12} & \alpha_{13}    \\
\alpha_{21}  & \alpha_{22}  & \alpha_{23}  \\
\alpha_{31}  & \alpha_{32} & \alpha_{33}
\end{array} \right).
$$
Comparing the coefficients, we deduce that
$R(e_{12}) = \alpha_{12}e_{12} + \alpha_{23}e_{23}$.

Similar arguments give us the following equalities:
$R(e_{21}) = \alpha_{21}e_{21} + \alpha_{32}e_{32}$,
$R(e_{23}) = \beta_{12}e_{12}  + \beta_{23}e_{23}$,
$R(e_{32}) = \beta_{21}e_{21}  + \beta_{32}e_{32}$,
$R(e_{13}) = \alpha_{13}e_{13}$,
$R(e_{31}) = \beta_{31}e_{31}$.

Consider $R(e_{13})$ and $R(e_{31})$. We have that
$$
0 = R(e_{13})R(e_{33})
  = (\alpha_{13}+1)R(e_{13}).
$$

Thus, if $R(e_{13})\neq0$, then $R(e_{13}) = -e_{13}$. Similarly, if
$R(e_{31})\neq0$, then $R(e_{31}) = -e_{31}$.
Since $\ker(R)$ and $\ker(R')$ are subalgebras in $M_3(F)$, we have that
$\alpha_{13}\neq\beta_{31}$.

Note that up to conjugation with transpose
it is sufficient to consider only one of these situations.
We assume that
$$
R(e_{13})=-e_{13},\quad R(e_{31})=0.
$$

Consider a subspace $M_1 = \Span\{e_{12},e_{23}\}$ spanned by $e_{12}$ and $e_{23}$.
As we know, $M_1$ is $R$-invariant.
Let $a\in M_1$ and
$R(a)=\gamma_{12}e_{12}+\gamma_{23}e_{23}$. We have
$$
\gamma_{23}e_{23}
 = R(a)R(e_{22})
 = R(\gamma_{12}e_{12}+a(e_{33}+e_{22}))
 = \gamma_{12}R(e_{12})+R(a).
$$

Therefore, $\gamma_{12}(R(e_{12})+e_{12}) = 0$ for all $a\in M_1$
and either $R(a)\in \Span\{e_{23}\}$ for all $a\in M_1$
or $R(e_{12})=-e_{12}$.

Similarly, if $M_2 = \Span\{e_{21},e_{32}\}$, then
either $R(a)\in \Span\{e_{32}\}$ for all $a\in M_2$
or $R(e_{21})=-e_{21}$.

If simultaneously $R(e_{12})=-e_{12}$ and $R(e_{21})=-e_{21}$,
then since $\ker(R')$ is a subalgebra in $M_3(F)$, $R(e_{11})=-e_{11}$, a~contradiction.

Suppose that $\alpha_{12}=\alpha_{21}=0$. Then
$$
\alpha_{23}\alpha_{32}e_{32}
 = R(e_{12})R(e_{21}) = R(e_{11}),
$$
a~contradiction. We have only two possibilities left:

Case 1a) $R(e_{21})=-e_{21}$, $R(e_{12})=\alpha_{23}e_{23}$, and
$R(e_{23})=\beta_{23}e_{23}$.

Case 1b) $R(e_{12})=-e_{12}$, $R(e_{21})=\alpha_{32}e_{32}$, and
$R(e_{32})=\beta_{32}e_{32}$.

Consider case 1a). Since $e_{21}$ and $e_{13}$ are in $\ker(R')$, then
so is $e_{23}=e_{21}e_{13}$. From
$$
0=R(e_{12})R(e_{33})=R(\alpha_{23}e_{23})=-\alpha_{23}e_{23},
$$
it follows that  $R(e_{12})=0$. Finally, since $e_{31}$ and $e_{12}$
are in $\ker(R)$, then $e_{32}=e_{31}e_{12}\in \ker(R)$. Summing up
the obtained equalities we have
$\{e_{13},e_{21},e_{23}\}=\ker(R')$ and $\{e_{12},e_{31},e_{23}\}\subset \ker(R)$.
It is a primitive RB-operator.

Consider case 1b). From
$0 = R(e_{33})R(e_{32})
   = (\beta_{32}+1)R(e_{32})$, it follows that
either $R(e_{32})=0$ or $R(e_{32})=-e_{32}$.

Suppose that $R(e_{32})=-e_{32}$.
For $R(e_{23})=\beta_{12}e_{12}+\beta_{23}e_{23}$, we have
$$
-\beta_{23}e_{22}
 = R(e_{23})R(e_{32})
 = R(\beta_{23}e_{22})
 = \beta_{23}e_{33}.
$$
Thus, $\beta_{23}=0$. Moreover, from
$0 = R(e_{22})R(e_{23}) = R(e_{23})$
it follows that $R(e_{23}) = 0$.
Since $\ker(R)$ is a subalgebra, $R(e_{21}) = 0$,
and we obtain a primitive RB-operator.

It remains to consider the case $R(e_{32}) = 0$. From
$$
0 = R(e_{23})R(e_{32})
  = R(\beta_{23}e_{22}+e_{22})
  = (\beta_{23}+1)e_{33},
$$
we obtain that $\beta_{23} = -1$.

Further, from the equation
$$
-\alpha_{32}e_{33}
 = R(e_{21})R(e_{23})
 = R(\alpha_{32}e_{33}+\beta_{12}e_{22})
 = \beta_{12}e_{33},
$$
we have $\alpha_{32} =-\beta_{12}$.
So, the RB-operator~$R$ satisfies
$e_{31},e_{32}\in \ker(R)$,
$e_{12},e_{13}\in\ker(R')$,
$R(e_{21})= -a e_{32}$, $R(e_{23}) = ae_{12} - e_{23}$.
If $a = 0$, we get a primitive RB-operator.
For $a\neq0$, we apply the conjugation with~$\Upsilon_a$,
\begin{equation}\label{Upsilon-a}
\begin{gathered}
\Upsilon_a(e_{ii}) = e_{ii},\ i=1,2,3,\quad
\Upsilon_a(e_{12}) = e_{12},\quad
\Upsilon_a(e_{21}) = e_{21},\\
\Upsilon_a(e_{13}) = a e_{13},\quad
\Upsilon_a(e_{23}) = a e_{23},\quad
\Upsilon_a(e_{31}) = (1/a) e_{31},\quad
\Upsilon_a(e_{32}) = (1/a) e_{32},
\end{gathered}
\end{equation}
and get the RB-operator 1-I).

\underline{Case 2)} $R(e_{11}) =  e_{22} + e_{33}$, $R(e_{22}) = e_{33}$, $R(e_{33}) = -e_{33}$.

By Lemma~1, the subalgebras
$M_1 = \Span\{e_{11},e_{22},e_{23},e_{32},e_{33}\}$,
$M_2 = \Span\{e_{12},e_{13}\}$,
$M_3 = \Span\{e_{21},e_{31}\}$
are invariant under the action of $R$,
and $R$ is splitting on each of $M_2$ and $M_3$.

Now, we use that $\ker(R')$ is homogeneous.
Since $e_{11}\not\in \ker(R)$ and $e_{11}\not\in \ker(R')$,
we have that
\begin{equation}\label{dimR=dimR'=2}
\dim (\ker(R)\cap M_2) + \dim (\ker(R)\cap M_3)
 = \dim (\ker(R')\cap M_2) + \dim (\ker(R')\cap M_3)
 = 2.
\end{equation}
We have two variants (up to conjugation with transpose).
{\bf Main variant}: $e_{12},e_{13}\in \ker(R')$ and $e_{21},e_{31}\in \ker(R)$.
We will consider it later.

{\bf Second variant}.
$e_{21},e_{13}\in \ker(R')$ and
$p = e_{12}+ye_{13}, q = e_{31}+ae_{21}\in \ker(R)$.
Since $\ker(R')$ is a subalgebra, $e_{23}\in\ker(R')$.
Let us multiply $p,q$, we will get again an element from $\ker(R)$,
$$
(e_{12}+ye_{13})(e_{31}+ae_{21})
 = (a+y)e_{11},
$$
so $y = -a$.
If $a = 0$, then $e_{12},e_{31}\in \ker(R)$ as well as $e_{32}\in \ker (R)$.
It is a primitive RB-operator.

Suppose that $a\neq0$. Thus,
\begin{equation}\label{r-element}
r = (e_{31}+ae_{21})(e_{12}-ae_{13})
 = ae_{22}-a^2e_{23}+e_{32}-ae_{33}\in \ker(R)
\end{equation}
and $R(e_{32}) = -a^2 e_{23} - 2a e_{33}$.
Then the RB-operator $P = \phi(\Upsilon_{1/a}^{-1}\Phi_{12}R\Phi_{12}\Upsilon_{1/a})$
is defined by Example~3 for $A_0 = D_3$, $A_- = U_3$
and $A_+ = \Span\{e_{21}-e_{23},e_{32}+e_{12},e_{11}-e_{33}+e_{31}-e_{13}\}$
with the action $P(e_{22}) = -e_{11} - e_{22} - e_{33}$,
$P(e_{11}) = -e_{11} - e_{33}$, $P(e_{33}) = 0$.

Define the automorphism $\varrho = \varrho(a,b,c)$ of $M_3(F)$ for nonzero $a,b\in F$ and any $c\in F$ as follows
\begin{equation}\label{2-III-Iso}
\begin{gathered}
\varrho(e_{11}) = e_{11}, \quad
\varrho(e_{12}) = \frac{1}{a}(e_{12} - bc e_{13}), \quad
\varrho(e_{13}) = b e_{13}, \\
\varrho(e_{21}) = a e_{21}, \quad
\varrho(e_{22}) = e_{22} - bc e_{23}, \quad
\varrho(e_{23}) = ab e_{23}, \quad
\varrho(e_{31}) = c e_{21} + \frac{1}{b}e_{31},\\
\varrho(e_{32}) = \frac{1}{a}\left( ce_{22}-b c^2 e_{23}+\frac{1}{b}e_{32}-c e_{33}\right),\quad
\varrho(e_{33}) = e_{33} + bc e_{23}.
\end{gathered}
\end{equation}
Conjugation with the automorphism~$\Phi_{12}\varrho\Phi_{12}$, where $a = b = c = 1$,
maps $P$ to a primitive RB-operator.

Let us return to the main variant.
By Lemma~5, the projection of $\Imm(R)$ on $e_{11}$ is zero.
So, the matrix algebra $N = \Span\{e_{22},e_{23},e_{32},e_{33}\}$ is $R$-invariant.
Now, we apply Corollary~1. In the case (M3) we get only primitive RB-operators.
In other cases, we have $R(e_{11}) = e_{22} + e_{33}$ and

2a) $e_{12},e_{13}\in \ker(R')$,
$e_{21},e_{31},e_{23},e_{22},e_{33}\in \ker(R)$,
$R(e_{32}) = e_{33} - e_{32}$ (by (M4${}^\prime$)),

2b) $e_{12},e_{13}\in \ker(R')$,
$e_{21},e_{31},e_{32},e_{22},e_{33}\in \ker(R)$,
$R(e_{23}) = e_{22} - e_{23}$ (by (M4${}^T$)),

2c) $e_{12},e_{13}\in \ker(R')$,
$e_{21},e_{31},e_{22},e_{33}\in \ker(R)$,
$R(e_{23}) = e_{22}- e_{23}$,
$R(e_{32}) = e_{33}- e_{32}$ (by (M5)),

2d) $e_{12},e_{13}\in \ker(R')$,
$e_{21},e_{31},e_{22},e_{33}\in \ker(R)$,
$R(e_{23}) = e_{33}- e_{23}$,
$R(e_{32}) = e_{22}- e_{32}$~(by (M5${}^T$)).

The following lemma shows that such trick is correct.

{\bf Lemma 6}.
Let $R_1$ and $R_2$ be RB-operators on $M_3(F)$ of weight~1 such that
$e_{12},e_{13}\in \ker(R_1'),\ker(R_2')$,
$e_{21},e_{31}\in \ker(R_1),\ker(R_2)$, and
$M = \Span\{e_{11}\}\oplus N$ is $R_1$- and $R_2$-invariant,
where $N = \Span\{e_{22},e_{23},e_{32},e_{33}\}$.
Suppose that the projections $R_1|_N$ and $R_2|_N$ of $R_1,R_2$ from $M$ on $N$
are conjugate with some automorphism of $N$.
Then there exists $\psi\in\Aut(M_3(F))$ such that the map
$Q = R_2 - \psi^{-1}R_1\psi$ is zero on~$M_2\oplus M_3$
and $\Imm(Q) \subset \Span\{e_{11}\}$.

{\sc Proof}.
Suppose that $R_1|_N$ and $R_2|_N$ as RB-operators on~$N$
are conjugate with an automorphism~$\xi$ of~$N$, i.e.,
$R_2|_N = \xi^{-1}R_1|_N\xi$.
Since each automorphism of $M_2(F)$ is inner, we have
$\xi(A) = T^{-1}AT$, $A\in N$, for some nondegenerate matrix $T\in N$.
Let us extend~$\xi$ from~$N$ on the entire algebra $M_3(F)$ as follows,
$$
\xi(X)
= \xi\left(\begin{pmatrix}
s & u \\
v & A
\end{pmatrix}\right)
 = \begin{pmatrix}
s & uT \\
T^{-1}v & T^{-1}AT
\end{pmatrix},
$$
where $s = x_{11}$, $u = (x_{12},x_{13})$, $v = (x_{21},x_{31})^T$.
It is an inner automorphism of $M_3(F)$ defined by the matrix
$\begin{pmatrix}
1 & 0 \\
0 & T
\end{pmatrix}$.

It is easy to check that
$Q = R_2 - \xi^{-1}R_1\xi$ is, maybe, nonzero only on~$N$
and $\Imm(Q) \subset \Span\{e_{11}\}$. Lemma is proved. \hfill $\square$

The RB-operator 2a) coincides with the RB-operator 5-I. The RB-operator 2b) is conjugate to the RB-operator 5-I with the help of $\Phi_{23}$.

The RB-operator 2c) is conjugate to the RB-operator 2-I) with the help of $\Phi_{23}\varrho\Phi_{23}$, where
$\varrho$ is defined by~\eqref{2-III-Iso} with $a = c = 1/b$.
Analogously, the RB-operator 2d) is conjugate to
the RB-operator 2-II).

\underline{Case 3)} $R(e_{11})= e_{22} + e_{33}$, $R(e_{22}) = -(e_{22}+e_{33})$, $R(e_{33})=-e_{33}$.

By Lemma~1, the subalgebras
$M_1 = \Span\{e_{11},e_{12},e_{21},e_{22},e_{33}\}$,
$M_2 = \Span\{e_{13},e_{23}\}$,
$M_3 = \Span\{e_{31},e_{32}\}$
are invariant under the action of $R$.
Let $N = \Span\{e_{11},e_{12},e_{21},e_{22}\}$.

Analogously to Case~2), we have two variants (up to conjugation with transpose).

{\bf Second variant}.
$e_{13},e_{32},e_{12}\in \ker(R')$ and
$p = e_{23}+ae_{13}, q = e_{31}-ae_{32}\in \ker(R)$.
In both cases $a = 0$ and $a\neq0$ we get only primitive RB-operators,
the proof is analogous to the one from Case 2).

{\bf Main variant}: $e_{13},e_{23}\in \ker(R')$ and $e_{31},e_{32}\in \ker(R)$.
Applying Statement~6, Corollary~1, and Lemma~6, we get the following RB-operators  defined by (M4${}^\prime$), (M4${}^T$), (M5), and (M5${}^T$) respectively (in (M3) we have only primitive ones),

3a) $R(e_{11}) = \alpha e_{33}$, $R(e_{12}) = 0$,
$R(e_{22}) = -\alpha e_{33}$, $R(e_{21}) = e_{22} - e_{21} + \beta e_{33}$;

3b) $R(e_{11}) = \alpha e_{33}$, $R(e_{21}) = 0$,
$R(e_{22}) = -\alpha e_{33}$, $R(e_{12}) = e_{11} - e_{12} + \beta e_{33}$;

3c) $R(e_{11}) = \alpha e_{33}$, $R(e_{12}) = e_{11}-e_{12}+\gamma e_{33}$, $R(e_{22}) = -\alpha e_{33}$, $R(e_{21}) = e_{22} {-} e_{21} {+} \beta e_{33}$;

3d) $R(e_{11}) = \alpha e_{33}$, $R(e_{12}) = e_{22}-e_{12}+\gamma e_{33}$, $R(e_{22}) = -\alpha e_{33}$, $R(e_{21}) = e_{11} {-} e_{21} {+} \beta e_{33}$.

\noindent
In all cases we have used that $R(e_{11}+e_{22}) = 0$ and the following easy fact,

{\bf Lemma 7}.
Let $R$ be an RB-operator of weight~1 on $M_3(F)$ and $R(e_{12}) = \alpha e_{12}+\beta e_{33}$, then $\beta = 0$.

{\sc Proof}. We are done by the equality
$\beta^2 e_{33}
 = R(e_{12})R(e_{12}) = 0$. \hfill $\square$

Note that in all cases 3a)--3d) the equality
$-\alpha^2e_{33} = R(e_{11})R(e_{22}) = 0$
implies $\alpha = 0$.

Let us consider the subcase 3a). From
\begin{equation}\label{3a,c}
\beta^2
 = R(e_{21})R(e_{21})|_{e_{33}}
 = R(e_{21})|_{e_{33}} = \beta,
\end{equation}
we have either $\beta = 0$ or $\beta = 1$.
If $\beta = 0$, then $R$ is splitting.
If $\beta = 1$, then $R$ is conjugate to the RB-operator 8-I) with the help of $\Phi_{13}\circ T$.

In the subcase 3b), we have an RB-operator which is conjugate to the RB-operator from 3a) with the help of $\Phi_{12}$.

Consider the subcase 3c). By~\eqref{3a,c}, we have $\beta\in\{0,1\}$. Further,
\begin{gather*}
\beta
 = R(e_{12}-e_{11})|_{e_{33}}
 = R(e_{12})R(e_{21})|_{e_{33}}
 = \beta\gamma
 = R(e_{21})R(e_{12})|_{e_{33}}
 = R(e_{21}-e_{22})|_{e_{33}} = \gamma.
\end{gather*}

If $\beta = \gamma = 0$, then $R$ is splitting. If $\beta = \gamma = 1$, $R$ is conjugate to the RB-operator 3-II) with the help of $\Phi_{13}\circ T\circ \varrho$, where $\varrho$ is defined by~\eqref{2-III-Iso} with $a=c=1/b$.

Analogously, in the subcase 3d) we get either splitting RB-operator or an RB-operator which is conjugate to the RB-operator 3-I) with the help of $\Phi_{13}\circ T\circ \varrho\circ\Phi_{13}$, where $\varrho$ is taken with $a=c=-1/b$.

\underline{Case 4)} $R(e_{11})=e_{22}+e_{33}$, $R(e_{22})=-e_{22}$, $R(e_{33})=0$.

As in Case~3), the subalgebras
$M_1 = \Span\{e_{11},e_{12},e_{21},e_{22},e_{33}\}$,
$M_2 = \Span\{e_{13}$, $e_{23}\}$,
$M_3 = \Span\{e_{31},e_{32}\}$
are invariant under the action of~$R$.

Analogously to Case~2), we have two variants (up to conjugation with transpose).

{\bf Second variant}.
$e_{13},e_{32},e_{12}\in \ker(R')$ and
$p = e_{23}+ae_{13}, q = e_{31}-ae_{32}\in \ker(R)$.
As above, we get only primitive RB-operators in both cases $a = 0$ and $a\neq0$.

{\bf Main variant}: $e_{13},e_{23}\in \ker(R')$ and $e_{31},e_{32}\in \ker(R)$.
Let us show that
\begin{equation}\label{Case4ImR}
\Imm(R) = \ker (R')\oplus \Span\{e_{33}\}
\end{equation}
(as vector spaces). Denote $A = \Imm(R)$. By Lemma~2,
$e_{11}A,Ae_{11},e_{22}A,Ae_{22}\subset \ker (R')$.
If
$a = \sum\limits_{i,j=1}^3\alpha_{ij}e_{ij}\in A$, then
$\alpha_{ij}e_{ij}\in \ker(R')$ for $i,j\in\{1,2\}$.
Subtracting them, we get that
$\alpha_{13}e_{13}$, $\alpha_{31}e_{31}$, $\alpha_{23}e_{23}$, $\alpha_{32}e_{32}\in \ker(R')$.
Thus, $A = \ker (R')\oplus \Span\{e_{33}\}$, otherwise $R$~is splitting.

Denote $N = \Span\{e_{11},e_{12},e_{21},e_{22}\}$.
As in Case~3), since $\ker(R')$ is homogeneous, we have the following cases:

a) $\ker(R')\cap N$ is one-dimensional with basis $e_{22}$;

b) $\ker(R')\cap N$ is two-dimensional with basis $e_{12},e_{22}$ or $e_{21},e_{22}$.

In the subcase a), we have
$R(e_{12}) = \alpha e_{22} + \beta e_{33}$ and
$R(e_{21}) = \gamma e_{22} + \delta e_{33}$.
From
$$
\alpha\gamma e_{22} + \beta\delta e_{33}
 = R(e_{21})R(e_{12})
 = R(e_{22}),
$$
we deduce that $\alpha\gamma = -1 $ and $\beta\delta = 0$.
From the projection of the equality
$R(e_{12})R(e_{21}) = R(\gamma e_{12} + \alpha e_{21} + e_{11})$
on the $e_{33}$-coordinate, we get $1 + \alpha\delta + \beta\gamma = \beta\delta$.
Thus, we obtain two RB-operators:

4a1) $e_{13},e_{23}\in \ker(R')$, $e_{31},e_{32}\in \ker(R)$,
$R(e_{12}) = a e_{22}$, $R(e_{21}) = (-1/a)(e_{22} + e_{33})$;

4a2) $e_{13},e_{23}\in \ker(R')$, $e_{31},e_{32}\in \ker(R)$,
$R(e_{21}) = a e_{22}$, $R(e_{12}) = (-1/a)(e_{22} + e_{33})$.

Conjugation with $\Phi_{13}\circ T\circ \Upsilon_{1/a}$
maps the RB-operator 4a1) to the RB-operator~$P$ with
$e_{12},e_{13},e_{22}\in\ker(P')$,
$e_{11},e_{21},e_{31}\in\ker(P)$,
$P(e_{23}) = e_{22}$, $P(e_{32}) = -e_{11} - e_{22}$,
$P(e_{33}) = e_{11} + e_{22}$.
Finally, $\psi^{-1}P\psi$ is the RB-operator 6-V).
Here $\psi = \Phi_{23}\varrho\Phi_{23}$, where~$\varrho$ is defined by~\eqref{2-III-Iso} with the parameters $a = c = -1/b$.

Conjugation with $\Phi_{13}\circ T\circ \Upsilon_a$,
maps the RB-operator 4a2) to the RB-operator~$P$ with
$e_{12},e_{13},e_{22}\in\ker(P')$,
$e_{11},e_{21},e_{31}\in\ker(P)$,
$P(e_{32}) = e_{22}$, $P(e_{23}) = -e_{11} - e_{22}$,
$P(e_{33}) = e_{11} + e_{22}$. Conjugation with $\varrho$ defined with $a = -c = -1/b$
maps $P$ to the RB-operator 6-IV).

Consider the subcase b). By Statement~6 and Corollary~1, we have the following subcases arisen from the cases (M5), (M5${}^T$), (M6), and (M6${}^T$) respectively,

\noindent4b1) $R(e_{11}) = \alpha e_{33}$, $R(e_{22}) = (1-\alpha)e_{33}$, $R(e_{12}) = e_{11}-e_{12}+\gamma e_{33}$, $R(e_{21}) = e_{22}-e_{21}+\delta e_{33}$;

\noindent4b2) $R(e_{11}) = \alpha e_{33}$, $R(e_{22}) = (1-\alpha)e_{33}$, $R(e_{21}) = e_{11}-e_{21}+\gamma e_{33}$, $R(e_{12}) = e_{22}-e_{12}+\delta e_{33}$;

\noindent4b3) $R(e_{11}) = -e_{11}+\alpha e_{33}$, $R(e_{22}) = e_{11}+(1-\alpha)e_{33}$, $R(e_{12}) = -e_{12}+\gamma e_{33}$, $R(e_{21}) = \delta e_{33}$;

\noindent4b4) $R(e_{11}) = -e_{11}+\alpha e_{33}$, $R(e_{22}) = e_{11}+(1-\alpha)e_{33}$, $R(e_{21}) = -e_{21}+\gamma e_{33}$, $R(e_{12}) = \delta e_{33}$.

We have used here that $R(e_{11}+e_{22}) = e_{33}$ for all variants.

Let us consider the case 4b1). From
\begin{gather}
\gamma e_{33}
 = R(e_{11}+e_{22})R(e_{12})
 = R(e_{12}-e_{12}+e_{11}) = \alpha e_{33}, \nonumber \\
 \delta e_{33}
 = R(e_{11}+e_{22})R(e_{21})
 = R(e_{21}-e_{21}+e_{22}) = (1-\alpha)e_{33},\nonumber \\
\alpha(1-\alpha)e_{33} \label{4b:Rf1Rf2}
 = R(e_{11})R(e_{22}) = 0,
\end{gather}
we get the RB-operators $R_1$ (when $\alpha=0$) and $R_2$ (when $\alpha=1$) satisfying
\begin{gather*}
R_1(e_{11}) = 0,\quad R_1(e_{22}) = e_{33},\quad
R_1(e_{12}) = e_{11}-e_{12},\quad R_1(e_{21}) = e_{22}-e_{21}+e_{33},\\
R_2(e_{11}) = e_{33},\quad R_2(e_{22}) = 0,\quad
R_2(e_{12}) = e_{11}-e_{12}+e_{33},\quad R_2(e_{21}) = e_{22}-e_{21}.
\end{gather*}
They are conjugate with the help of $\Phi_{12}$. The conjugation of $R_1$ with $\Phi_{13}\circ T\circ \varrho$, where $\varrho$ is defined by~\eqref{2-III-Iso} with $a = -c = 1/b$, gives the RB-operator 4-I).

In the case 4b2) we analogously get the RB-operator 4-II).

Consider the case 4b3). By
$$
\gamma e_{33}
 = R(e_{11}+e_{22})R(e_{12})
 = R(e_{12}-e_{12}) = 0, \quad
 \delta^2 e_{33}
 = R(e_{21})R(e_{21}) = 0,
$$
we get only primitive RB-operators. We deal analogously with the subcase 4b4).

\underline{Case 5)}
$R(e_{11}) = e_{22} + e_{33}$, $R(e_{22}) = R(e_{33}) = 0$.

We follow the strategy from Case~2).
Since $e_{11}\not\in \ker(R)$ and $e_{11}\not\in\ker(R')$, we get~\eqref{dimR=dimR'=2}.

We have $\ker(R)$-homogeneity.
So, we have two variants (up to conjugation with transpose).

{\bf Second variant}.
$e_{21},e_{13}\in\ker(R)$ and so, $e_{23}\in \ker(R)$. Also,
$p = e_{12}-ae_{13}, q = e_{31}+ae_{21}\in \ker(R')$.
As above, we get only primitive RB-operators in both cases $a = 0$ and $a\neq0$.

{\bf Main variant}: $e_{12},e_{13}\in \ker(R')$ and $e_{21},e_{31}\in \ker(R)$.
By Lemma~5, the projection of $\Imm(R)$ on $e_{11}$ is zero.
So, the matrix algebra $N = \Span\{e_{22},e_{23},e_{32},e_{33}\}$ is $R$-invariant.
Suppose that $e_{23},e_{32}\in\ker(R)$, then we get the RB-operator~5-I).

Then we have three subcases:

5a) $e_{23}\in\ker(R)$;

5b) $e_{32}\in\ker(R)$;

5c) $e_{23},e_{32}\not \in \ker(R)$.

In 5a--5c), we apply Corollary~1 (joint with Lemma~6) to get the following RB-operators (in (M3) we get only primitive ones)

5-II. $e_{12},e_{13}\in \ker(R')$,
$e_{21},e_{31},e_{23}\in \ker(R)$,
$R(e_{32}) = e_{33} - e_{32}$ (by (M4${}^\prime$)),

5-III. $e_{12},e_{13}\in \ker(R')$,
$e_{21},e_{31},e_{32}\in \ker(R)$,
$R(e_{23}) = e_{22} - e_{23}$ (by (M4${}^T$)),

5-IV. $e_{12},e_{13}\in \ker(R')$,
$e_{21},e_{31}\in \ker(R)$,
$R(e_{23}) = e_{22}- e_{23}$,
$R(e_{32}) = e_{33}- e_{32}$ ((M5)),

5-V. $e_{12},e_{13}\in \ker(R')$,
$e_{21},e_{31}\in \ker(R)$,
$R(e_{23}) = e_{33}- e_{23}$,
$R(e_{32}) = e_{22}- e_{32}$~((M5${}^T$))

\noindent respectively.

Conjugation of the RB-operator 5-IV) with the automorphism
$\Phi_{23}\varrho\Phi_{23}$ coincides with the RB-operator from the case 2-I).
Here $\varrho$ is defined by~\eqref{2-III-Iso} with $a = c = 1/b$.
Also, conjugation of the RB-operator 2-II) with the same $\varrho$
gives the RB-operator 5-V).

The RB-operator 5-III) is conjugate
to the RB-operator 5-II) with the help of $\Phi_{23}$.

\underline{Case 6)} $R(e_{11}) = e_{22}$, $R(e_{22}) = R(e_{33}) = 0$.

In this case  $R(1)=e_{22}$ and by Lemma~1 we obtain the following equations:
\begin{gather*}
R(e_{12})=\alpha_{12}e_{12}+\alpha_{32}e_{32},\quad
R(e_{21})=\alpha_{21}e_{21}+\alpha_{23}e_{23},\\
R(e_{23})=\beta_{21}e_{21}+\beta_{23}e_{23},\quad
R(e_{32})=\beta_{12}e_{12}+\beta_{32}e_{32},\\
R(e_{13})=\alpha_{11}e_{11}+\alpha_{13}e_{13}+\alpha_{22}e_{22}+\alpha_{31}e_{31}+\alpha_{33}e_{33},\\
R(e_{31})=\beta_{11}e_{11}+\beta_{13}e_{13}+\beta_{22}e_{22}+\beta_{31}e_{31}+\beta_{33}e_{33}.
\end{gather*}
From
$$
\alpha_{22}e_{22}=R(e_{11})R(e_{13})=R(\alpha_{11}
e_{11}+(\alpha_{13}+1)e_{13}),
$$
we conclude that
\begin{equation}\label{us1}
\text{if}\  \alpha_{13}=-1,\  \text{then}\  \alpha_{11}=\alpha_{22}.
\end{equation}

Similarly,
\begin{equation}\label{us2}
\text{if}\ \beta_{31}=-1,\  \text{then}\  \beta_{11}=\beta_{22}.
\end{equation}

Consider a subalgebra
$M_1 = \Span\{e_{11},e_{22},e_{33},e_{13},e_{31}\}$.
We proved that $M_1$ is $R$-invariant.
It is easy to see that
$M_1 = \Span\{e_{22}\}\oplus N$ (as vector spaces) for
$N = \Span\{e_{11},e_{33},e_{13},e_{31}\}\cong M_2(F)$.

Let $M_2 = \Span\{e_{21},e_{23}\}$ and $M_3 = \Span\{e_{12},e_{32}\}$.
Then $M_2$ and $M_3$ are $R$-invariant subspaces in $M_3(F)$.
From~(7) it follows that restrictions $R_i$, $i=2,3$, on $M_i$ satisfy
$R_i^2 + R_i=0$. By Lemma~2, $e_{33}(R(e_{32})+e_{32})\in \ker(R)$.
Thus, $(\beta_{32}+1)R(e_{32})=0$ and
\begin{equation}\label{us3}
\text{if}\  R(e_{32})\neq 0,\  \text{then}\ \beta_{32}=-1.
\end{equation}
Similarly,
\begin{equation}\label{us4}
R(e_{23})=0\ \text{or}\  \beta_{23}=-1.
\end{equation}

Also, from
$$
0 = R(e_{21})R(e_{11})=(\alpha_{21}+1)R(e_{21}),\quad
0 = R(e_{11})R(e_{12})=(\alpha_{12}+1)R(e_{12})
$$
we obtain that
\begin{equation}\label{us5}
R(e_{21})=0\ \text{or}\  \alpha_{21}=-1,\quad
R(e_{12})=0\ \text{or}\  \alpha_{12}=-1.
\end{equation}

Let $Pr_N$ be the projection from $M_1$ onto $N$ with
$\ker(Pr_N) = \Span\{e_{22}\}$. By Statement 6 the composition
$\bar{R} = Pr_N\circ R$ is the Rota-Baxter operator on $N$.
Note that
$\bar{R}(e_{11})=\bar{R}(e_{33})=0$.
Now we apply Corollary~1 and obtain that $\bar{R}$ may be conjugate to (M3), (M4) or (M5).

As in Lemma 6, we can find an RB-operator $Q$ satisfying the conditions\\
(i) $Q$ is conjugate to $R$,\\
(ii) $M_1$, $M_2$, and $M_3$ are $Q$-invariant,\\
(iii) the projection $\bar{Q}$ is equal to (M3), (M4) or (M5),\\
(iv) $Q(e_{22})=0$.

Note that in all cases (M3), (M4), and (M5), $Pr_N(e_{11})=Pr_N(e_{33})=0$. Therefore, $Q(e_{11})=\delta_1 e_{22}$ and $Q(e_{33})=\delta_2 e_{22}$.

Let us show that we may assume that\\
(v) $Q(e_{11}) = e_{22}$ and $Q(e_{33}) = 0$.

From
$$
\delta_1\delta_2 e_{22} = Q(e_{11})Q(e_{33}) = 0,
$$
it follows that $\delta_1=0$ or $\delta_2=0$.
Since $R$ and $Q$ are conjugate with automoprhism or antiautomorphism of $M_3(F)$,  $\tr(Q(1)) = \tr(R(1)) = 1$. Therefore, $\delta_1=1$ and $\delta_2=0$ or vice versa.

It is straightforward to check that in Cases (M3) and (M5), up to conjugation with $T\circ\Phi_{13}$  and $\Phi_{13}$ respectively, we can assume that $\delta_1=1$, $\delta_2=0$, i.e., the condition (v) is fulfilled.

Suppose that the restriction $\bar{Q}$ corresponds to (M4) and $\delta_1=0$, $\delta_2=1$. Then we have
\begin{gather*}
Q(e_{11})=0,\quad
Q(e_{22})=0, \quad
Q(e_{33})=e_{22},\quad
Q(e_{13})=0,\quad
Q(e_{31})=e_{11}-e_{31}+\gamma e_{22}
\end{gather*}
for some $\gamma\in F$.
Consider the RB-operator $Q_1=(T\circ\Phi_{13})^{-1}Q(T\circ\Phi_{13})$. We have
\begin{gather*}
Q_1(e_{11})=e_{22},\quad
Q_1(e_{22})=0, \quad
Q_1(e_{33})=0,\quad
Q_1(e_{13})=0,\quad
Q_1(e_{31})=e_{33}-e_{31}+\gamma e_{22}.
\end{gather*}
For $Q_1$, we may apply \eqref{us2} and obtain $\gamma=0$. Define an automorphism $\psi$ of $M_3(F)$ as
\begin{gather*}
\psi(e_{11})=e_{11}-e_{13},\
\psi(e_{22})=e_{22},\
\psi(e_{33})=e_{33}+e_{13},\
\psi(e_{12})=e_{12},\
\psi(e_{21})=e_{21}-e_{23},\\
\psi(e_{13})=e_{13},\quad
\psi(e_{31})=e_{11}-e_{33}+e_{31}-e_{13},\quad
\psi(e_{23})=e_{23},\quad
\psi(e_{32})=e_{32}+e_{12}.
\end{gather*}
\noindent
It is easy to check that the RB-operator $Q_2=\psi^{-1}Q_1\psi$ satisfies all conditions (i)--(v).

Consider the following subcases:\\
\emph{Subcase 6a}: $R(e_{21})=0$;\\
\emph{Subcase 6b}: $R(e_{21})\neq 0$ and $R(e_{23})\neq 0$;\\
\emph{Subcase 6c}:  $R(e_{21})\neq 0$ and $R(e_{23})= 0$.

\emph{Subcase 6a}) Suppose that $R(e_{21})=0$. Since $e_{11},e_{31}\notin \ker(R)$, then
$e_{12},e_{32}\notin \ker(R)$.
By Lemma~3, $\ker(R)$ is homogeneous, therefore the restriction
$R_3$ on $M_3$ is non-degenerate.
From $R_3^2+R_3=0$ we derive $R_3=-\id$ in this case.
Thus, $R(e_{12})=-e_{12}$ and $R(e_{32})=-e_{32}$.
Finally, for any $\beta\in F$ we have
$$
(-\beta_{21}-\beta \beta_{23})e_{22}
 = R(e_{23})R(e_{12}+\beta e_{32})
 = R((\beta_{21}+\beta \beta_{23})e_{22})=0.
$$
Therefore, $\beta_{21}=\beta_{23}=0$ and $R(e_{23})=0$.

Suppose that the restriction $\bar{R}$ is equal to (M3) from Corollary~1. In this case,
$$
R(e_{13})=\gamma_1 e_{22},\ \ R(e_{31})=-e_{31}+\gamma_2 e_{22}.
$$
By Lemma 7, $\gamma_1=\gamma_2=0$ and $R$ is equal to the RB-operator 6b).

Suppose that the restriction $\bar{R}$ is equal to (M5). Then we have:
\begin{gather*}
R(e_{11})=e_{22},\quad
R(e_{22})=0, \quad
R(e_{33})=0,\quad
R(e_{12})=-e_{12},\quad
R(e_{21})=0,\\
R(e_{13})=e_{11}+\gamma_1e_{22}-e_{13},\\
R(e_{31})=-e_{31}+e_{33}+\gamma_2 e_{22},\quad
R(e_{23})=0,\quad
R(e_{32})=-e_{32}.
\end{gather*}

We can now use \eqref{us1} and \eqref{us2} to get $\gamma_1=1$ and $\gamma_2=0$.

Define an automorphism $\xi$ of $M_3(F)$ as
\begin{gather*}
\xi(e_{11})=e_{22},\
\xi(e_{22})=e_{11}+e_{13},\
\xi(e_{33})=e_{33}-e_{13},\
\xi(e_{12})=e_{21}+e_{23},\
\xi(e_{21})=e_{12},\\
\xi(e_{13})=e_{23},\quad
\xi(e_{31})=e_{32}-e_{12},\quad
\xi(e_{23})=e_{13},\quad
\xi(e_{32})=e_{33}-e_{13}-e_{11}+e_{31}.
\end{gather*}
\noindent The conjugation of $R$ with $T\circ\xi$ gives us the RB-operator 4-I.
If the restriction $\bar{R}$ corresponds to (M4), then similar reasons give us that
\begin{gather*}
R(e_{11})=e_{22},\
R(e_{22})=R(e_{33})=0,\
R(e_{12})=-e_{12},\\
R(e_{21})=0,\
R(e_{13})=0,\
R(e_{31})=e_{11}+e_{22}-e_{31},\
R(e_{23})=0,\
R(e_{32})=-e_{32}.
\end{gather*}
This operator is conjugate to the RB-operator 6-IV with the help of $\Phi_{12}\circ T$.

\emph{Subcase 6b}) Now suppose that $R(e_{21})\neq 0$ and $R(e_{23})\neq 0$.
Then by Lemma~3 the restriction $R_2$ is non-degenerate.
Therefore, $R(e_{21})=-e_{21}$ and $R(e_{23})=-e_{23}$.

If $R(e_{12})\neq 0$, then by \eqref{us5}, $R(e_{12})=-e_{12}+\alpha_{32}e_{32}$.
But from
\begin{equation}\label{6b:e21-e_12}
e_{22} = R(e_{21})R(e_{12}) = R(-e_{22}) = 0
\end{equation}
we get a contradiction. Thus, $R(e_{12})=0$.
By the same reasons, inequality $R(e_{32})\neq 0$ contradicts with
$$
e_{22} = R(e_{23})R(e_{32}) = R(-e_{22}) = 0.
$$

As in Subcase 6a, we consider the restriction of $R$ on $N$ and obtain RB-operators corresponding to the cases (M3)--(M5).
Using similar arguments as in 6a we obtain that if the restriction $\bar{R}$ corresponds to (M3), we obtain operator 6c). If the restriction corresponds to (M5), then $R$ lies in the orbit of 4-II and if the restriction corresponds to (M4), then $R$ is conjugate to the RB-operator 6-V.

\emph{Subcase 6c}) The only remaining question is what happens if
$R(e_{21})\neq 0$ and $R(e_{23}) {=} 0$.

Suppose that $\bar{R}$ corresponds to (M3). Then, by Lemma 7,  $R(e_{31})=-e_{31}$ and \linebreak
$R(e_{13})=0$.

Consider $R(e_{21})$. By \eqref{us5}, $R(e_{21})=-e_{21}+\alpha_{23}e_{23}$. From
$$
-\alpha_{23}e_{21}=R(e_{21})R(e_{31})=R(\alpha_{23}e_{21})
$$
we obtain that $\alpha_{23}=0$ and $R(e_{21}) = -e_{21}$.

We will consider two subcases: $R(e_{32})=0$ and $R(e_{32})\neq 0$. If $R(e_{32})=0$, then, since $\ker(R)$ is a subalgebra in $M_3(F)$, $R(e_{12})=0$. Conjugation with $T$ gives us the RB-operator 6-VI.

If $R(e_{32})\neq 0$, then from \eqref{us3} we have $R(e_{32})=-e_{32}+\beta_{12}e_{12}$. From
$$
-\beta_{12}e_{32}=R(e_{31})R(e_{32})=R(\beta_{12}e_{32})
$$
we deduce that $R(e_{32})=-e_{32}$. By~\eqref{6b:e21-e_12}, we get $R(e_{12}) = 0$ and therefore $R$ is a primitive RB operator.

Suppose that the restriction $\bar{R}$ corresponds to (M5). Then $R$ satisfies
$R(e_{13})=-e_{13}+ e_{11}+\beta e_{22}.$
From \eqref{us1} we deduce that
$R(e_{13})=-e_{13}+e_{11}+ e_{22}$.
The same reasons give us
$R(e_{31})=-e_{31}+e_{33}$.

Recall that $R(e_{21})\neq 0$ implies $R(e_{21})=-e_{21}+\alpha_{23}e_{23}$ by~\eqref{us5}. From
$$
- e_{21}+e_{23}=R(e_{21})R(e_{13})=R(-e_{23}+ e_{21}
-e_{23}+e_{23})= R(e_{21}),
$$
we obtain that
$R(e_{21})=-e_{21}+e_{23}$.

Since $R(e_{23})=0$ and $R(e_{13})\neq 0$, then $R(e_{12})\neq 0$.
Thus, by \eqref{us5}, $R(e_{12})=-e_{12}+\alpha_{32}e_{32}$. From
$$
(1+\alpha_{32})e_{22}
 = R(e_{21})R(e_{12})=R(-e_{22}-e_{22}+e_{22}) = 0,
$$
we deduce that $\alpha_{32}=-1$ and $R(e_{12})=-e_{12}- e_{32}$.
Finally, since $(R^2+R)(e_{12})=0$ we obtain that  $R(e_{32})=0$.
Thus,  we obtain the RB-operator 6-I.

Now suppose that the restriction $\bar{R}$ corresponds to (M4).
We have $R(e_{13})=0$,
$R(e_{31})=-e_{31}+e_{11}+e_{22}$,
$R(e_{21})=-e_{21}+\alpha_{23}e_{23}$, and $R(e_{23})=0$.

From
\begin{multline*}
- e_{21}-\alpha_{23}e_{21}
 = R(e_{21})R(e_{31}) \\
 = R(\alpha_{23}e_{21}+ e_{21})
 = (\alpha_{23}+1)R(e_{21})=(\alpha_{23}+1)(-e_{21}+\alpha_{23}e_{23}),
\end{multline*}
we obtain that
$\alpha_{23}(\alpha_{23}+1) = 0$.
Therefore, $\alpha_{23} = -1$ or $\alpha_{23} = 0$.

Consider $e_{12}$.
If $R(e_{12})\neq 0$, then
$R(e_{12})=-e_{12}+\alpha_{32}e_{32}$ by~\eqref{us5}. From
\begin{multline*}
e_{11}-\alpha_{23}e_{13}-\alpha_{32}e_{31}+\alpha_{23}\alpha_{32}e_{33}=R(e_{12})R(e_{21}) \\
 = R(-e_{11}+\alpha_{32}e_{31}-e_{11}+\alpha_{23}e_{13}+e_{11})=-e_{22}+\alpha_{32}R(e_{31}),
\end{multline*}
we obtain $\alpha_{32}\neq 0$. The last inequality holds if and only if $R(e_{32}) = 0$,
otherwise the restriction $R_2$ is non-degenerate and consequently $R_2 =-\id$.
Since $\ker(R)$ is a~subalgebra, $e_{13}e_{32}=e_{12}\in \ker(R)$, a~contradiction.

Therefore, $R(e_{12}) = 0$.
If $R(e_{32}) = 0$, then we obtain  the RB-operators 6-II (when $\alpha_{23}=0$) and 6-III (when $\alpha_{23}=-1$).

Suppose that $R(e_{32})\neq 0$. Then
$R(e_{32}) = -e_{32}+\beta_{12}e_{12}$. Consider
\begin{multline*}
e_{31}-\alpha_{23}e_{33}-\beta_{12}e_{11}+\beta_{12}\alpha_{23}e_{13}=R(e_{32})R(e_{21}) \\
 = R(-e_{31}+\beta_{12}e_{11}-e_{31}+\alpha_{23}e_{33}+e_{31})=-R(e_{31})+\beta_{12}e_{22}.
\end{multline*}

So,
$R(e_{31})=\beta_{12}(e_{11}+e_{22})-e_{31}+\alpha_{23}e_{33}-\beta_{12}\alpha_{23}e_{13}$.
It holds if and only if $\beta_{12} = 1$ and $\alpha_{23}=0$.
Therefore,  $e_{21},\  e_{32}- e_{12}\in\ker(R')$ and
we obtain $e_{21}(e_{32}- e_{12}) = - e_{22}\in\ker(R')$, a contradiction.

\underline{Case 7)}
$R(e_{11}) = e_{22}$, $R(e_{22}) = 0$, $R(e_{33}) = -e_{33}$.

We have $R(1) = e_{22}-e_{33}$ and by Lemma~1 (as in Case~1) we deduce
\begin{equation}\label{Case7Variables}
\begin{gathered}
R(e_{12}) = \alpha_{12}e_{12}+\alpha_{31}e_{31},\quad
R(e_{31}) = \beta_{12}e_{12}+\beta_{31}e_{31}, \quad
R(e_{23}) = \gamma_{23}e_{23}, \\
R(e_{21}) = \alpha_{21}e_{21}+\alpha_{13}e_{13}, \quad
R(e_{13}) = \beta_{21}e_{21}+\beta_{13}e_{13}, \quad
R(e_{32}) = \gamma_{32}e_{32}.
\end{gathered}
\end{equation}

Moreover, by Lemma~1, $\gamma_{23},\gamma_{32}\in\{0,-1\}$.
Since $\ker(R)$ and $\ker(R')$ are subalgebras in $M_3(F)$,
then $\gamma_{23}\neq\gamma_{32}$.
Up to conjugation with transpose, we may assume that
$$
R(e_{23})=-e_{23},\quad R(e_{32})=0.
$$

From the equation
$$
0 = R(e_{12})R(e_{22})=(\alpha_{12}+1)R(e_{12}),
$$
it follows that either $\alpha_{12}=-1$ or $R(e_{12})=0$.
Similarly, either $\alpha_{21}=-1$ or $R(e_{21})=0$.
The case $\alpha_{12} = \alpha_{21} = 0$
does not hold, otherwise $e_{11}\in \ker(R)$.
The case $\alpha_{12} = \alpha_{21} = -1$
does not hold, otherwise $e_{11}\in \Imm(R)$,
a contradiction to~\eqref{Case7Variables}.
So, we have two subcases:

7a) $\alpha_{12}=-1$, $R(e_{21})=0$.
Since $\ker(R)$ is a subalgebra in $M_3(F)$,
then $R(e_{31})=0$. From
$$
-\alpha_{31}e_{31}=R(e_{33})R(e_{12})=R(\alpha_{31}e_{31})=0,
$$
it follows that $R(e_{12})=-e_{12}$ and since $\ker(R')$
is a~subalgebra, then $R(e_{13})=-e_{13}$.
We obtain that $R$ is a primitive RB-operator.

7b) $\alpha_{21}=-1$, $R(e_{12})=0$.
Since $R(e_{21})\neq 0$, then
$0 = R(e_{22})R(e_{13})=R(\beta_{21}e_{21})$ implies $\beta_{21}=0$.
From
$0=R(e_{11})R(e_{13})=(\beta_{13}+1)R(e_{13})$, it follows that either
$R(e_{13})=-e_{13}$ or $R(e_{13})=0$.

From
$$
\beta_{12}e_{12}=R(e_{31})R(e_{11})=(\beta_{31}+1)R(e_{31}),
$$
it follows that either $R(e_{31})=-e_{31}$ or $\beta_{31}=0$.
If $\beta_{31} = 0$, then
$$
-\beta_{12}e_{11}
 = R(e_{31})R(e_{21})
 = R(\beta_{12}e_{11}+\alpha_{13}e_{33})
 = \beta_{12}e_{22}-\alpha_{13}e_{33}
$$
lead us to $\beta_{12} = \alpha_{13} = 0$.
Thus, either $R(e_{31})=-e_{31}$ or $R(e_{31}) = 0$, $R(e_{21}) = -e_{21}$.

If $R(e_{13}) = 0$, then, since $\ker(R)$ is a~subalgebra, $R(e_{31})=-e_{31}$.
So, $e_{21} = e_{23}e_{31}\in\ker(R')$.
Applying the conjugation with $\Phi_{23}$, we get that $R$ is a primitive RB-operator.

If $R(e_{13}) = -e_{13}$, then $R(e_{31}) = 0$ and so, $R(e_{21}) = -e_{21}$.
Applying the conjugation with $\Phi_{12}$, we get that
$R$~is again a~primitive RB-operator.

\underline{Case 8)}
$R(e_{11}) = -e_{11}$, $R(e_{22}) = R(e_{33}) = 0$.

We follow the strategy from Case~2). Since $e_{11}\not\in \ker(R)$
and $e_{22},e_{33}\not\in\ker(R')$, we again get~\eqref{dimR=dimR'=2}.

We have $\ker(R)$-homogeneity. So, we have
{\bf Main variant}: $e_{12},e_{13}\in \ker(R')$ and $e_{21},e_{31}\in \ker(R)$.
We will consider it later.

{\bf Second variant}.
$e_{21},e_{13},e_{23}\in \ker(R)$ and
$p = e_{12}-ae_{13}, q = e_{31}+ae_{21}\in \ker(R')$.
In both cases when $a = 0$ or $a\neq0$, we get primitive RB-operators
(the proof is analogous to the one from Case~2)
which are splitting in Case 8).

Let us return to the main variant.
Suppose that $e_{23},e_{32}\in\ker(R)$, then $R$ is splitting.
Then we have three subcases:

a) $e_{23}\in \ker(R)$,

b) $e_{32}\in\ker(R)$,

c) $e_{23},e_{32}\not \in \ker(R)$.

Let us consider the subcase a).
By Lemma~1, $Fe_{11}\oplus N$ is $R$-invariant, where
$N = \Span\{e_{22},e_{23},e_{32},e_{33}\}$.
By Statement~6 and Corollary~1, we may assume that
$R(e_{32})
 = \left(\begin{matrix}
\gamma_{11} & 0 & 0 \\
0 & \gamma_{22} & 0 \\
0 & -1 & 0
\end{matrix}\right)$, where $\gamma_{22}\in\{0,1\}$.
From the equality
$R(e_{32})R(e_{32}) = \gamma_{22}R(e_{32})$,
we get $\gamma_{11}(\gamma_{11}-\gamma_{22}) = 0$.
If $\gamma_{11} = 0$, we get a splitting RB-operator.
So, $\gamma_{11} = \gamma_{22} = 1$
gives the RB-operator~8-I.

In b), analogously to a), we get the RB-operator~$R$ satisfying
$e_{12},e_{13}\in \ker(R')$, $e_{21},e_{31},e_{32}\in \ker(R)$,
$R(e_{23}) = e_{11}+e_{33} - e_{23}$,
it is conjugate to the RB-operator 8-I) with the help of~$\Phi_{23}$.

In c), denote
$R(e_{23}) = (\delta_{ij})$ and
$R(e_{32}) = (\gamma_{ij})$.
By~Statement~6, the projection $R|_N$ of $R$ on the subalgebra~$N$ is an RB-operator.
Applying Corollary~1, we have either
$R(e_{23}) = \delta_{11}e_{11} + e_{22} - e_{23}$ and
$R(e_{32}) = \gamma_{11}e_{11} - e_{32} + e_{33}$ or
$R(e_{23}) = \delta_{11}e_{11} - e_{23} + e_{33}$ and
$R(e_{32}) = \gamma_{11}e_{11} + e_{22} - e_{32}$.
We consider the first variant, the second one is similar.
From
\begin{gather*}
\gamma_{11}\delta_{11}e_{11}+e_{22}-e_{23}
 = R(e_{23})R(e_{32})
 = R(e_{23}), \\
\gamma_{11}\delta_{11}e_{11}+e_{33}-e_{32}
 = R(e_{32})R(e_{23})
 = R(e_{32}),
\end{gather*}
we get that either $\gamma_{11} = \delta_{11} = 0$ and so, $R$ is splitting, or
$\gamma_{11} = \delta_{11} = 1$. So, we obtain the RB-operator~$R$
satisfying $e_{12},e_{13}\in \ker(R')$, $e_{21},e_{31}\in \ker(R)$,
$R(e_{23}) = e_{11} + e_{22} - e_{23}$,
$R(e_{32}) = e_{11} + e_{33} - e_{32}$.
The RB-operator~$R$ is conjugate to the RB-operator~3-II) with the help of
$\varrho$ is defined by~\eqref{2-III-Iso} with $a=c=1/b$.

Analogously, for the second variant we get either splitting RB-operator or
the RB-operator~$R$ satisfying $e_{12},e_{13}\in \ker(R')$,
$e_{21},e_{31}\in \ker(R)$,
$R(e_{32}) = e_{11} + e_{22} - e_{32}$,
$R(e_{23}) = e_{11} + e_{33} - e_{23}$
which is conjugate to the RB-operator~3-I).

\underline{Case 9)}
$R(e_{11}) = e_{22}$, $R(e_{22}) = -e_{22}$, $R(e_{33}) = -e_{33}$.

{\bf Lemma 8}.
In Case 9), there are no non-splitting RB-operators.

Denote $A = \Imm(R) = R(M_3(F))$.
As in Case 4), we get that either $A = \ker(R+\id)$ and then $R$ is splitting
or $A = \ker (R+\id)\oplus\Span\{e_{11}\}$ (as vector spaces).

Consider the second variant. Suppose that $R(y) = e_{11}$ for some $y\in M_3(F)$. So,
\begin{equation}
0 = R(y)R(1)
  = R(e_{11} + y(-e_{33}) + y)
  = e_{11} + e_{22} - R(ye_{33}). \label{9:UnitOnpreImagef1}
\end{equation}
From~\eqref{9:UnitOnpreImagef1}, we have $R(ye_{33}) = e_{11}+e_{22}$.
Let $y = \alpha e_{13}+\beta e_{23} + \gamma e_{33} + y'$
for $y'\in\Span\{e_{ij}\mid i=1,2,3,\,j=1,2\}$. Thus,
$$
R(\alpha e_{13}+\beta e_{23} + \gamma e_{33} - e_{11}) = e_{11},
$$
and we may assume that
$y = \alpha e_{13}+\beta e_{23} + \gamma e_{33} - e_{11}$.
Further,
\begin{multline*}
e_{11}
 = R(y)R(y)
 = R(e_{11}y + ye_{11} + y^2) \\
 = R(-e_{11}+\alpha e_{13} -e_{11} + \alpha\gamma e_{13}+\beta\gamma e_{23}+\gamma^2 e_{33}-\alpha e_{13}+e_{11} ) \\
 = R((\gamma-1)e_{11}+\gamma(\alpha e_{13}+\beta e_{23} + \gamma e_{33} - e_{11}))
 = (\gamma-1)e_{22} + \gamma e_{11},
\end{multline*}
which leads us to $\gamma = 1$.

We may rewrite the formula $R(y) = e_{11}$ as
$$
1 = R(\alpha e_{13}+\beta e_{23}-e_{11}-e_{22})
  = R(\alpha e_{13}+\beta e_{23}).
$$
Define $z = \alpha e_{13}+\beta e_{23}$.
On the one hand,
$$
1 = R(z)R(z)
  = R(R(z)z+zR(z)+z^2) = 2 + R(z^2)
$$
and so $R(z^2) = -1$. On the other hand, $z^2 = 0$.
We have a contradiction. \hfill $\square$

We have considered all cases of the action $R$ on $D_3(F)$.
Theorem is proved. \hfill $\square$

{\bf Corollary 2}.
All RB-operators obtained in Theorem~3 lie in different orbits
under the action of the operator $\phi$ from Statement 1 and conjugation with automorphisms of $M_3(F)$ and transpose.

{\sc Proof}.
Let us denote by $X*$ the set of cases $Xa,Xb,Xc$ for $X\in\{1,2,3,4,6,7\}$
and $Xa,Xb$ for $X = 5$.

The Jordan form of $R$ as a linear map as well
as rank of $R(1)$ have to be preserved under the action of $\Aut(M_3(F))$ and transpose.
So, we may compute the 6-tuple
\begin{equation}\label{6-tuple}
(\dim(\ker(R)),\dim(\ker(R^2)),
\dim(\ker(R')),\dim(\ker(R'^2)),
\rank(R(1)),\rank(R'(1)))
\end{equation}
for each case and compare them up to the action of~$\phi$.

Immediately, we get that the cases
\begin{center}
$1*)$, 1-I), 5-I), $7*)$, 8-I)
\end{center}
lie in their own orbits.
Indeed, the cases $1*)$ and 1-I) are unique with the property $R^2(R+\id)^2\neq 0$,
and they lie in different orbits, since their minimal polynomials
$m_1 = x^3(x+1)$ and $m_{1-I} = x^3(x+1)^2$ do not coincide.
Only for the cases $7*)$, we have $(\rank (R(1)),\rank (R'(1))) = (2,2)$.
The case 8-I) is unique case with $\dim(\ker R^2) = \dim(\ker R) = 5$.
The case 5-I) is the only case satisfying the conditions $\dim(\ker(R)) = 6$ and $(\rank (R(1)),\rank (R'(1)) = (2,3)$.

Let us check that all listed in Theorem~3 primitive RB-operators lie in different orbits.

{\bf Lemma 9}.
The cases $Xa$, $Xb$ and $Xc$ for $X\in\{1,2,3,4,6,7,8,9\}$
and $5a$, $5b$ lie in pairwise different orbits.

{\sc Proof}.
Analyzing 6-tuples~\eqref{6-tuple}, we obtain that primitive RB-operators of different types
lie in different orbits.

Consider two different RB-operators $O$ and $P$ of the same type~$X$.
Suppose that they lie in the same orbit, so
either $\psi^{-1}O\psi = P$ or $\psi^{-1}O\psi = T\circ P\circ T$ for a~$\psi\in\Aut(M_3(F))$.
Consider the first case.
Note that $\psi$ preserves both kernels and their powers.
Moreover, $\psi$~preserves the radicals of the kernels and the one-dimensional annihilators of such radicals, i.e.,
$$
\psi(e_{13}) = se_{13},\quad
\psi(e_{31}) = te_{31}
$$
for some $s,t\in F$. So, $\psi(e_{11}) = st e_{11}$ and $\psi(e_{33}) = st e_{33}$.
Since the image of an idempotent under the action of an automorphism has to be an idempotent, $st = 1$.
Thus, $\psi(e_{11}) = e_{11}$, $\psi(e_{33}) = e_{33}$, and $\psi(e_{22}) = e_{22}$,
since $\psi$ preserves the identity matrix. So, we get a contradiction.

In the second case, we analogously get that
$$
\psi(e_{13}) = se_{31},\quad
\psi(e_{31}) = te_{13}.
$$
Again, $st=1$ and $\psi(e_{11}) = e_{33}$, $\psi(e_{33}) = e_{11}$, and $\psi(e_{22}) = e_{22}$.
It means that the action of $O$ and $P$ on the subalgebra of diagonal matrices
should coincide up to the action of $\Phi_{13}$. We get a contradiction for the all cases $Xa,Xb,Xc$.
\hfill $\square$

Let us continue on the separation of cases by 6-tuples~\eqref{6-tuple}.
Up to $\phi$ we have the following orbits,

a) $2*)$, 2-I), 2-II) for the 6-tuple $(4,5,4,4,2,3)$,

b) $3*)$, 3-I), 3-II) for the 6-tuple $(4,4,4,5,1,2)$,

c) $4*)$, 4-I), 4-II), 6-I) for the 6-tuple $(4,5,4,4,1,3)$,

d) 6-II), 6-III), 6-VI) for the 6-tuple $(6,7,2,2,1,3)$,

e) $5*)$, 5-II) for the 6-tuple $(5,6,3,3,2,3)$,

f) $6*)$, 6-IV), 6-V) for the 6-tuple $(5,6,3,3,1,3)$.

Let us show how the analysis of left and right annihilators of the kernels
helps to separate cases.
If RB-operators $P$ and $Q$ lie in the same orbit, then both their
kernels should be isomorphic or anti-isomorphic.
In particular, pairs of the dimensions of left and right annihilators
for each of $\ker(R)$ and $\ker(R+\id)$ should be pairwise equal.

For example, in the case 6-IV) we have
$\ker(R_{6-IV})=\Span\{e_{11},e_{33},e_{32},e_{21},e_{31}\}$ and
$\Ann_l(\ker(R_{6-IV})) = \Ann_r(\ker(R_{6-IV})) = (0)$.

In the case 6-V) we have
$\ker(R_{6-V})=\Span\{e_{11},e_{33},e_{23},e_{21},e_{31}\}$ and
$\Ann_l(\ker(R_{6-V}))$\\
$=(0)$ but
$\Ann_r(\ker(R_{6-V}))=\Span\{e_{23},e_{21}\}$.
It means that 6-IV and 6-V lie in different orbits.

By the same argument applied for $\ker(R+\id)$, we separate cases 2-I) and 2-II),
4-I) and 4-II) respectively.

a) Further, the cases $2*)$ do not lie in same orbit with neither 2-I) nor 2-II), since
$\ker(R_{2*}^2)$ has two-dimensional semisimple part, while
$\ker(R^2)$ has a three-dimensional semisimple part for $R$ in the cases 2-I) and 2-II).

b) The semisimple part of $\ker(R)$ is one-dimensional in the cases $3*)$
and it is two-dimensional in the cases 3-I) and 3-II).
We will show below that the RB-operators from the cases 3-I) and 3-II)
lie in different orbits.

c) The case 6-I) lies in its own orbit since it is
the only variant from c) satisfying the condition $\ker(R)\cong M_2(F)$.
The cases $4*)$ do not lie in the same orbit with neither of 4-I), 4-II).
Indeed, $\ker(R_{4*})$ has one-dimensional semisimple part,
when the semisimple part of $\ker(R)$ in the cases 4-I), 4-II) is two-dimensional.

d) The subalgebra $\ker(R+\id)$ has trivial product only in the case 6-VI).
Further, $\ker(R+\id)^2 = \ker(R+\id)$ only in the case 6-II).

e) The cases $5*)$ and 5-II) do not lie in the same orbit.
The algebra $\Imm(R)$ in the case $5*)$ has one-dimensional semisimple part,
when the semisimple part of $\Imm(R)$ in the case 5-II) is two-dimensional.

f) The cases $6*)$ do not lie in the same orbit with neither 6-IV) nor 6-V),
since $\ker(R'_{6*})$ is nilpotent when both $\ker(R'_{6-IV})$ and $\ker(R'_{6-V})$ have one-dimensional semisimple part.

Finally, let us show that the RB-operators $P$ and $O$ taken from the cases 3-I and 3-II respectively lie in different orbits. Assume there exists $\psi\in\Aut(M_3(F))$ such that
$\psi^{-1}O\psi = P$ or $\psi^{-1}O\psi = T\circ P\circ T$.
Note that in both cases $\psi(P(1)) = O(1)$. So, $\psi(e_{11}) = e_{11}$.
Consider the first case. We have
\begin{gather*}\allowdisplaybreaks
\ker(P+\id)   = \Span\{e_{11},e_{12},e_{13},e_{23}\},\quad
\ker(O+\id)   = \Span\{e_{11},e_{12},e_{13},e_{32}\},\\ \allowdisplaybreaks
\ker(P+\id)^2 = \ker(P+\id) \oplus \Span\{e_{22}\},\quad
\ker(O+\id)^2 = \ker(O+\id) \oplus \Span\{e_{22}\},\\
\ker(P) = \Span\{e_{21},e_{31},e_{22}+e_{33},e_{22}+e_{32}\},\
\ker(O) = \Span\{e_{21},e_{31},e_{22}+e_{33},e_{22}+e_{23}\}.
\end{gather*}
Further, $\psi(e_{13}) = ae_{12}$ for some nonzero $a\in F$, since
$\psi$ has to map the centralizer of $\rad(\ker(P+\id))$ onto the centralizer of $\rad(\ker(O+\id))$.
As $\psi$ maps $\rad(\ker P)$ onto $\rad(\ker O)$, we get
$$
\psi(e_{21}) = be_{21} + ce_{31},\quad
\psi(e_{31}) = de_{21} + fe_{31}.
$$
Further, $\psi(e_{33}) = \psi(e_{31})\psi(e_{13}) = e_{22} + af e_{32}$, so
$\psi(e_{22}) = \psi(e_{21})\psi(e_{12}) = e_{33} - af e_{32}$.
The following equalities
$$
(P+\id)(e_{22})
 = - e_{11}
 \neq 1
 = \psi^{-1}(1)
 = \psi^{-1}(O+\id)(e_{33} - af e_{32})
 = \psi^{-1}(O+\id)\psi(e_{22})
$$
imply a contradiction.

Now consider the second case when $\psi^{-1}O\psi = T\circ P\circ T$.
Then analogously, we get $\psi(e_{31}) = a e_{12}$ for some nonzero $a\in F$, and
$$
\psi(e_{12}) = be_{21} + ce_{31}, \quad
\psi(e_{13}) = de_{21} + fe_{31}.
$$
Thus,
$\psi(e_{33}) =  \psi(e_{31})\psi(e_{13}) = e_{11}$, a contradiction with $\psi(e_{11}) = e_{11}$.
\hfill $\square$

{\bf Remark 3}.
One can derive from Theorem~3 the classification of all non-splitting RB-operators of nonzero weight on the 5-dimensional semisimple associative algebra $A = Fe\oplus M_2(F)$, here $e^2 = e(\neq0)$. Indeed, given an RB-operator $R$ of weight one on $A$, we may extend its action on the entire algebra $M_3(F)$ by Example~2. More detailed, we embed $A$ into $M_3(F)$ as follows: $\psi(e) = e_{11}$, $\psi(e_{ij}) = e_{i+1\, j+1}$ for $e_{ij}\in M_2(F)$. Then we put $e_{12},e_{13}\in\ker(R+\id)$ and $e_{21},e_{31}\in\ker(R)$. If one starts with a~non-splitting RB-operator $R$ on~$A$, then its extension $R$ on $M_3(F)$ is again a~non-splitting RB-operator. So, $R$ up to $\phi$ and up to conjugation with an automorphism of $M_3(F)$ and transpose is one of the RB-operators from Theorem~3. On the other hand, all RB-operators from Theorem~3 except the cases 1-I), 6-I), 6-II), 6-III) are exactly mentioned above extensions of RB-operators on~$A$.

\section*{Acknowledgements}

M. Goncharov was supported by Russian Scientific Fond (project N 19-11-00039).

V. Gubarev was supported by the Program of fundamental scientific researches of the Siberian Branch of Russian Academy of Sciences, I.1.1, project 0314-2019-0001.

\noindent Maxim Goncharov \\
Vsevolod Gubarev \\
Novosibirsk State University \\
Pirogova str. 2, 630090 Novosibirsk, Russia \\
Sobolev Institute of Mathematics \\
Acad. Koptyug ave. 4, 630090 Novosibirsk, Russia \\
e-mail: gme@math.nsc.ru, wsewolod89@gmail.com

\end{document}